\input amstex
\documentstyle{amsppt}
\document
\magnification=1200
\NoBlackBoxes
\nologo
\pageheight{18cm}


\bigskip

\centerline{\bf REAL MULTIPLICATION} 

\smallskip

\centerline{\bf AND NONCOMMUTATIVE GEOMETRY}

\smallskip 

\centerline{\it (ein Alterstraum)}

\medskip

\centerline{\bf Yu.~I.~Manin}

\medskip

\centerline{\it Max--Planck--Institut f\"ur Mathematik, Bonn}

\bigskip

{\bf Abstract.} Classical theory of Complex Multiplication (CM)
shows that all abelian extensions of a complex quadratic field $K$
are generated by the values of appropriate modular
functions at the points of finite order of elliptic
curves whose endomorphism rings are orders in $K$.
For real quadratic fields, a similar description is not known.
However, the relevant (still unproved) case of Stark conjectures ([St1])
strongly suggests that such a description must exist. 
In this paper we propose to use 
two--dimensional quantum tori corresponding to real quadratic
irrationalities as a replacement of elliptic curves
with complex multiplication. We discuss some basic
constructions of the theory of
quantum tori from the perspective of this Real Multiplication (RM)
research project.

\bigskip

\hfill{$ \overset , \to \epsilon\pi\epsilon\grave{\iota}\ \kappa\alpha
\grave{\iota}\ \tau\grave{\alpha}\
\gamma\nu\acute{\omega}\rho\iota\mu\alpha\ 
\overset , \to o\lambda\acute{\iota}\gamma o\iota\varsigma\
\gamma\nu\acute{\omega}\rho\iota\mu\alpha\ 
\overset , \to \epsilon\sigma\tau\iota\nu$}

\smallskip

\hfill{\it ... for even subjects that are known are known only to a few.}

\smallskip

\hfill{\it Aristotle, Poetics IX, 1451b}

\medskip

\centerline{\bf Preface}

\medskip

This paper can be read as a union of three largely independent 
parts.

\smallskip

Section 1 is dedicated to a general problem of noncommutative
geometry Connes style: what are
morphisms between noncommutative spaces considered
as ``spectra'' of associative rings (perhaps, with an 
additional structure)? One natural suggestion is to
define morphisms as isomorphism classes
of biprojective bimodules as in Morita theory. 
Slightly extending
Rieffel's Morita classification of
 two--dimensional
quantum tori, I present a 
description of the resulting category (Theorem 1.7.1)
in terms of what can be called
``period pseudolattices'' (sec. 1.1), by analogy with period
lattices of elliptic curves.
(In the context of operator algebras, requiring
quite sophisticated modification of basic notions,
A.~Connes calls such morphisms ``correspondences'',
cf. [Co1], p. 526, and [Jo2]).

\smallskip

Section 2 contains some results on the values and residues
of zeta functions of arithmetical progressions
in real quadratic fields, in the
spirit of earlier work of E.~Hecke, continued by G.~Herglotz
and D.~Zagier. Our calculations are strongly motivated
by H.~M.~Stark's conjectures ([St1], [St2]) proposing very special
generators of abelian extensions of such fields.  

\smallskip

Section 3 is a contribution to the theory of quantum theta functions
(see [Ma3]). It gives a partial answer to the question of A.~Schwarz
([Sch2]) about the relationship between quantum thetas
and representations of quantum tori. The main Theorem 3.7 of
this section generalizes a seminal calculation of F.~Boca in [Bo2].

\smallskip

I collected these disjoint results under one roof
because I feel that they form pieces of a general
picture, which could be called Real Multiplication
of two--dimensional quantum tori, by analogy with
the classical Complex Multiplication of elliptic curves
(Kronecker's {\it Jugendtraum}).

\smallskip

From this perspective, Section 1 outlays basics of the
(noncommutative) geometry of Real Multiplication, Section 2 presents
certain known or conjectural arithmetical facts
in the light of this geometric picture, whereas
Section 3 provides
elements of function theory. 
\smallskip

Unfortunately, the relations between these parts that I can establish
are too sparse yet. An important test for such a theory
would be a proof of Stark's conjectures
for real quadratic fields. If this plan succeeds,
a more ambitious project could address
Real Multiplication of multidimensional
quantum tori, as an analog and an extension of Shimura--Taniyama
multidimensional CM theory.  

\smallskip

I tried to facilitate reading this paper for potential
readers with varied backgrounds by providing
many definitions and introductory explanations,
so that large parts of this paper can be read as a review.
In particular, the introductory Section 0 explains
rudiments of the classical Complex Multiplication
theory which serves as a guide for the whole
enterprise. It discusses as well a 
possibilty of including this theory
in the context
of noncommutative geometry. For additional connections,
see [Ma4].

\medskip

{\it Acknowledgement.} The crystallization of this project
owes much to Matilde Marcolli and our collaboration
[MaMar]. Victor Nistor consulted me about the proof of
Lemma 1.4.2. Florian Boca's paper [Bo2] and correspondence
with him were crucial for recognizing the connection
between quantum thetas and modules over quantum tori.
Sasha Rosenberg's insights about morphisms between 
noncommutative spaces developed in [Ro2] helped me to
overcome a difficult psychological barrier.
Last but not least, I appreciate the proposal
of Friedrich Hirzebruch to translate {\it  Alterstraum}
in the title as ``midlife crisis''.

\bigskip


\centerline{\bf \S 0. Introduction: Lattices, elliptic curves,}

\smallskip

\centerline{\bf and Complex Multiplication}

\medskip

{\bf 0.0. An overview.} Let $K$ be a field of algebraic numbers of one of the three types: $\bold{Q}$, a complex quadratic extension
of $\bold{Q}$, or a real quadratic extension of $\bold{Q}.$
Consider the following classical problem: describe the maximal
abelian extension $K^{ab}$ of $K$. Of course, the Galois group
of such an extension is known  for arbitrary agebraic number fields
$K$: it is the id\`ele class group of $K$ modulo its
connected component. However, explicit generators of $K^{ab}$
and the action of the Galois group on them generally remain a mystery,
with exception of two classical cases described below.
 
\smallskip

According to the Kronecker--Weber theorem (KW),
$\bold{Q}^{ab}$ is generated by roots of unity, i.e.
by the points of finite order of the multiplicative group
$\bold{G}_m$ considered as an algebraic group over $\bold{Q}.$
For $K$ imaginary quadratic, the multiplicative group
should be replaced by the elliptic curve $E_K$ whose
$\bold{C}$--points are $\bold{C}/O_K$, $O_K$ being the ring of integers in $K$.
To get $K^{ab}$, one must adjoin to $K$ the values of a power
of the Weierstrass function at points of finite order of
$E_K$, and the value of the absolute invariant of $E_K$.
To see that points of finite order generate
an abelian extension, one observes that the action
of the Galois group on them must commute with the action
of algebraic endomorphisms furnished by the power maps
$x\mapsto x^m$ in the KW case, resp. the complex multiplication (CM)
maps written additively on the universal covering of
$E_K$ as $x\mapsto ax,\,a\in O_K.$ The commutant
of this action suitably completed in profinite topology is
abelian, and essentially coincides
with the completion of the action itself. The 
universal id\`elic description of the Galois group
together with reduction modulo $p$ arguments furnish the rest.

\smallskip

Elliptic curves have a rich analytic theory. Curves admitting
a complex multiplication form a subfamily of all elliptic curves.
The latter can be parametrized by their period lattices $\Lambda$ i.e. 
discrete images of the injective
homomorphisms
$j:\,\bold{Z}^2\to\bold{C}$ modulo a natural equivalence relation.
The moduli space of them is $PGL(2,\bold{Z})\setminus (H^+\cup H^-)$,
$H^{\pm}$ being the upper/lower  halfplanes respectively. 
The curves isogeneous to $E_K$ live over orbits
of points $\bold{P}^1(K)$. The multiplicative group
also appears in this family as the ``degenerate elliptic curve''
over the cusp, that is the orbit 
$PGL(2,\bold{Z})\setminus \bold{P}^1(\bold{Q})$, so that in principle
the geometry of the CM and KW cases can be unified. 
 
\smallskip

The cusp corresponds to the very degenerate lattice: $j$ acquires
a cyclic kernel. There is an intermediate case of degeneration,
invisible in algebraic geometry,
where $j$ is still injective, but its image is not discrete.
The relevant modular orbit is  $PGL(2,\bold{Z})\setminus (\bold{P}^1(\bold{R}) 
\setminus \bold{P}^1(\bold{Q}))$, it contains
orbits of $\bold{P}^1(K)$ for real quadratic $K$, but they could not be used
in the same way as CM points of the modular curve because
of lack of the analog of elliptic curves over this stratum
of the moduli space. Hopefully, quantum tori might serve
as a substitute.

\smallskip

This introductory section is dedicated to some details of
the CM picture and its possible extension to the RM case.

\medskip

{\bf 0.1. Category of lattices $\Cal{L}$.} By definition, {\it a lattice}
(of rank 2) is a triple $(\Lambda ,V,j)$, where 
$\Lambda$ is a free abelian group of rank two, $V$ is
an one--dimensional complex space, and $j:\,\Lambda\to V$ is an injective homomorphism
with discrete image, hence compact quotient.

\smallskip

When no confusion is likely, we will refer to $(\Lambda ,V,j)$
simply as $\Lambda$.

\smallskip

{\it A morphism}
of lattices $(\Lambda^{\prime},V^{\prime},j^{\prime})\to (\Lambda ,V,j)$ is a commutative diagram
$$
\CD 
\Lambda^{\prime} @>j^{\prime}>> V^{\prime}\\
@V\varphi VV  @VV\psi V  \\
\Lambda @>>j> V  \\
\endCD
\eqno(0.1)
$$
in which $\varphi$ is a group homomorphism, and
$\psi$ is a $\bold{C}$--linear map. Clearly, $\varphi$
is uniquely determined by $\psi$, and vice versa.
Choosing a basis $(\lambda_1,\lambda_2)$ in $\Lambda$, 
taking $j(\lambda_2)$ as the base vector of $V$  we see that
in any isomorphism class of lattices 
one can find a representative given by $j:\,\bold{Z}^2\to\bold{C}$
such that $j(0,1)=1$, $j(1,0):=\tau$ is a number in
$\bold{C}\setminus \bold{R}.$ Changing the sign of $\lambda_1$
if needed we can arrange $\tau$ to lie in the upper
half--plane $H$.

\smallskip

Let us denote this lattice $\Lambda_{\tau}.$
Then any non--zero morphism $\Lambda_{\tau^{\prime}} \to \Lambda_{\tau}$ is represented
by a non--degenerate matrix 
$$
g=\left(\matrix a&b\\c&d\endmatrix\right)\in M\,(2,\bold{Z})
$$
such that
$$
 \tau^{\prime} =\frac{a\tau +b}{c\tau +d} \, .
\eqno(0.2)
$$ 
This $g$ is obtained by writing  $\varphi$ in (0.1)
as the right multiplication of a row by a matrix;
the respective $\psi$ is the multiplication by $(c\tau +d)^{-1}.$

\smallskip

Clearly, (0.2) is an isomorphism, iff $g\in GL(2,\bold{Z})$. 
Thus the moduli space of (isomorphism classes of) lattices
is
$$
PGL(2,\bold{Z})\setminus (\bold{P}^1(\bold{C})\setminus \bold{P}^1(\bold{R}))
=PSL(2,\bold{Z})\setminus H .
\eqno(0.3)
$$
Endomorphisms of a lattice $(\Lambda ,V,j)$ form a ring,
with componentwise addition of $(\phi, \psi)$ and
composition as multiplication. It contains 
$\bold{Z}$ and comes together with its embedding in $\bold{C}$:
$$
\roman{End}\,\Lambda = \{a\in\bold{C}\,|\,aj(\Lambda )\subset j(\Lambda )\}.
$$ 
\proclaim{\quad 0.1.1. Lemma} (a) $\roman{End}\,\Lambda \ne \bold{Z}$
iff there exists a complex quadratic subfield $K$ of $\bold{C}$ 
such that $\Lambda$ is isomorphic to a lattice contained in $K.$

\smallskip

(b) If this is the case, denote by $O_K$ the ring of integers of
$K$. There exists a unique integer $f\ge 1$ (conductor) such that
$\roman{End}\,\Lambda =\bold{Z}+fO_K=:R_f,$ and $\Lambda$ is
a projective module of rank 1 over $R_f$. Every $K$, $f$
and a projective module over $R_f$ come from a lattice.

\smallskip

(c) If lattices $\Lambda$ and $\Lambda^{\prime}$ have the same
$K$ and $f$, they are isomorphic if and only if
their classes in the Picard group $\roman{Pic}\,R_f$
coincide.
\endproclaim

Automorphisms of a lattice generally form a group $\bold{Z_2}$
($\psi$ is multiplication by $\pm 1$.) However, integers
of two imaginary quadratic fields obtained by adjoining
to $\bold{Q}$ a primitive root of unity of degree 4 (resp. 6)
furnish examples of lattices with automorphism group
of order 4 (resp. 6). Only these two fields produce
lattices with such extra symmetries.

\medskip

{\bf 0.2. Category of elliptic curves $\Cal{E}$.} For any lattice 
$(\Lambda ,V,j)$
the quotient space $V/j(\Lambda )$ is an one--dimensional complex torus
which has a canonical structure of (the set of complex points of) an algebraic curve $E_{\Lambda}$ of genus 1 with base point $0$.
Such curves form a category $\Cal{E}$ (morphisms should respect
base points). 

\medskip

{\bf 0.3. The functor $P:\,\Cal{E} \to \Cal{L}.$} 
Let $E$ be an elliptic curve. The functor ``period lattice''
$P$ is defined on objects by the following prescription:
$P(E)=(\Lambda_E,V_E,j_E)$ where $V_E$ = the tangent space to $E$ at the base point,
considered as its Lie algebra, $\Lambda_E$ the kernel of the
of the exponential map $V_E\to E(\bold{C})$, and $j_E$ its
canonical embedding. On morphisms, $\psi$ is the induced tangent
map and $\varphi$ its restriction to the period lattices.

\medskip

\proclaim{\quad 0.3.1. Theorem} $P$ is an equivalence of categories.
\endproclaim    

\smallskip

This simple result is crucial for the theory of complex multiplication.

\medskip

{\bf 0.4. Abelian extensions of complex quadratic fields.}
 Let now $K$ be a complex
quadratic extension of $\bold{Q}$. Choose and fix an
embedding $K\to\bold{C}.$ Denote by $O_K$
the ring of integers of $K$. 

\smallskip

There are three related
but somewhat different ways to describe the maximal abelian extension $K^{ab}$
of $K$.

\medskip

(A) {\it Approach via elliptic curves.}

\smallskip

Here one starts with a single elliptic curve $E_K$ associated
to the lattice $O_K\subset \bold{C}$. It turns out
that its minimal definition field containing $K$
is generated by the
value of its absolute invariant $J(E_K)$, and is the maximal
unramified extension of $K$. One can also give a beautiful
description of the total set of conjugates of $J(E_K)$
and the action of the Galois group on this set.
Namely, any lattice whose endomorphism ring is precisely
$O_K$, is represented by an ideal in $O_K$, and two lattices
are isomorphic iff they lie in the same class.
Absolute invariants of the respective elliptic curves
are conjugate to each other, and the action
of the Galois group is induced by a geometric
twisting operation producing from a curve an isogenous
curve.

\smallskip

The remaining part of the $K^{ab}$ is generated
by the values at points of finite order of $E$ of 
a special function $t$. In Weierstrass notation, it is
$t=\wp (z,O_K)^{u}$ where $u$ is the order
of the automorphism group of $E$, so that our algebraic
numbers can be described as the values of a transcendental function
$$
\left[ \frac{1}{z^2}+
\sum_{\lambda\in O_K\setminus \{0\}}\left( \frac{1}{(z+\lambda )^2} -
 \frac{1}{\lambda^2}\right)\right]^u,
\eqno(0.4)
$$
at $z\in K.$
In geometric terms,
$t$ is an appropriate coordinate on the projective line
$$ 
\bold{P}^1=E_K/O_K^*,
\eqno(0.5)
$$
to which the points of finite order are mapped.

\medskip

(B) {\it Approach via modular curves.} 

\smallskip

In another approach,
one considers an extension of $K$ generated
by roots of unity and absolute invariants of all elliptic curves
admitting complex multiplication by an order
in $O_K$. Not all of $K^{ab}$ is generated
in this way, it remains to produce an additional
infinite extension with a Galois group of period two,
but in a sense the most essential part of $K^{ab}$
is obtained in this way.

\smallskip

This approach stresses the geometry
and arithmetic of {\it the moduli space (stack)} of elliptic
curves rather than that of elliptic curves themselves.
This space has special points which can be characterized
as fixed points of certain correspondences,
and fields of definition of these points are of primary
interest.

\smallskip

For a brief introduction to both approaches,
see [Se] and [Ste].
 
\medskip

(C) {\it Approach via Stark's numbers.} 

\smallskip

The general conjectures
due to H.~M.~Stark provide (hypothetical) generators
of abelian extensions which are values of zetas
(or their derivatives, or the Taylor coefficients next to the residue)
similar to (0.4). In the CM case these conjectures are
proved in [St2], by reducing them to the more classical and geometrical
forms of the theory sketched in (A), (B).
No independent arguments are known.

\smallskip 

To provide the basis
for comparison with the RM case, we will briefly
describe these numbers.

\smallskip

Let $(\Lambda ,\bold{C},j)$ be a lattice in $\bold{C}$,
$\lambda_0 \in \Lambda\otimes \bold{Q}.$
Put
$$
\zeta (\Lambda , \lambda_0, s):= \sum_{\lambda\in \Lambda}
\frac{1}{|j(\lambda_0 +\lambda )|^{2s}}
\eqno(0.6)
$$
where $j$ is extended by $\bold{Q}$--linearity. These series admit
meromorphic continuation and may have a pole
of the first order at $s=1$ and zero of the first order
at $s=0$. Similar behavior is exhibited in the real case.
Zeta-functions for two isomorphic lattices
differ by a factor $A^s$ where $A$ is a positive real number.
In the CM case we will restrict the choice
of $\Lambda$ in the isomorphism class by considering only
lattices with $j(\Lambda )\subset K.$ Then $A$ can be modulus squared of
any number in $K$.

\smallskip

We have the following simple lemma.

\smallskip

\proclaim{\quad 0.4.1. Lemma} (a) Assume that $F(s)$ 
vanishes at $s=0$. Then for any $A>0$
$$
\left.\frac{d}{ds}F(s)\right|_{s=0} =
 \left.\frac{d}{ds}(A^sF(s))\right|_{s=0} .
$$
In particular, 
$$
S_0(\Lambda,\lambda_0):=e^{\zeta^{\prime}(\Lambda , \lambda_0, 0)}
\eqno(0.7)
$$
is an invariant of the isomorphism class of $(\Lambda , \lambda_0)$.

\smallskip
(b) Assume that 
$$
F(s)= \frac{r}{s-1} +v +O(s-1),\ r\ne 0
$$
near $s=1$. Then the similar formula holds for $A^sF(s)$,
with the ratio $v/r$  replaced by $v/r+\roman{log}\,A$.
In particular, the following coset
$$
S_1(\Lambda,\lambda_0 ):= e^{v/r}\,\roman{mod}\,N_{K/\bold{Q}}(K^*)\in 
\bold{C}^*/N_{K/\bold{Q}}(K^*)
\eqno(0.8)
$$
is an invariant of the isomorphism class of $(\Lambda , \lambda_0)$.
Here $v,r$ are calculated via $\zeta (\Lambda , \lambda_0, s)$
for any representative of this class
satisfying $j(\Lambda )\subset K$, 
\endproclaim

\smallskip

The essence of Stark's conjectures consists in the prediction that
invariants  of the type $S_0(\Lambda,\lambda_0)$ are algebraic
units in appropriate abelian extensions of $K$, and that the action
of the Frobenius elements of the Galois group upon them can be 
explicitly described.

\smallskip

Stark's proof in the CM case is based upon a direct calculation of these
invariants, which in turn reduces 
to the second Kronecker limit formula. A version involving
$S_1(\Lambda,\lambda_0)$ might be more feasible from the computational
viewpoint, the two versions being essentially equivalent 
thanks to the classical functional
equations. These calculations show that Stark's numbers a priori defined as
values of some transcendental functions admit
an algebraic geometric interpretation demonstrating their
arithmetical nature.

\medskip

{\bf 0.5. Real quadratic fields.} In this paper we propose 
some constructions parallel to (A) -- (C) above, for the case
of real quadratic fields.

\smallskip

(A) {\it Geometry of real multiplication.} 

\smallskip

Replacing lattices 
by pseudolattices and elliptic curves by quantum tori,
we develop in in \S 1 the geometric framework parallel to that of
0.1--0.3 above.

\smallskip

(B) {\it Geometry of noncommutative modular curves.}

\smallskip

The space $PGL(2,\bold{Z})\setminus (\bold{P}^1(\bold{R}) 
\setminus \bold{P}^1(\bold{Q}))$ as an invisible stratum
of the classical modular curve was studied in [MaMar].
In particular, it was shown that its $K$--theory
can be written in terms of modular symbols,
and that classical modular forms of weight two and their Mellin transforms
are represented by interesting densities on this stratum.

\smallskip

For the purposes of real multiplication, however,
more relevant might be noncommutative spaces which
represent the orbits $PGL(2,\bold{Z})\setminus \bold{P}^1(K)$
corresponding to the  individual real quadratic $K$.
It seems that the remarkable paper by Bost and Connes
([BoCo]) and its extensions [ArLR], [HaL], furnish the right language
to describe the arithmetic phenomena that interest us.
Provisionally, [BoCo] appears to describe the KW case (cusp)
from the noncommutative viewpoint. However, a satisfactory
generalization of [BoCo] to more general number fields
is not developed as yet (cf. however [HaL], [ArLR], [Coh1]).

\smallskip

We expect that other noncommutative spaces, besides
quantum tori and modular curves, must play an 
essential role in the future
theory. In particular, the projective line (0.5)
might be replaced by the crossed product
of the algebra of functions on $K$
and its automorphism group of the type
$x\mapsto ax+b$ where $a\in O_K^*,\,b\in O_K.$
This looks even closer to the spaces studied in
[BoCo] and [HaL].

\medskip

(C) {\it Stark's numbers.} 

\smallskip

For real quadratic fields $K$, 
one should consider series
of the type (0.6), in which $|a|^{-2s}=(a\overline{a})^{-s}$
is replaced by $N_{K/\bold{Q}}(a)^{-2s}$ and furnished
with a slight additional twist: the typical
term of (0.6) is multiplied by the sign of the conjugate
of $j(l_0+l)$, cf. \S 2 below for more details.
More important is the following complication:
the action of an infinite cyclic group of units
makes each term in (0.6) to repeat infinitely often,
so to make sense of the whole expression one should only
sum over cosets modulo the relevant group. 

\smallskip

We show in \S 2 that an adaptation
of Hecke's calculations leads to formulas for Stark's numbers
which are compatible with the general picture
of ``passing to the quantum limit'': cf. 1.8 below.

\smallskip

Our hope, based upon this calculation, is that
an appropriate algebraic geometric refinement of (A)
and (B) will lead to a proof of Stark's conjectures,
in the same way as it worked in the CM case. 

\smallskip

Again, a natural question arises: can one see 
in noncommutative
geometry Stark's numbers
of the cyclotomic (KW) case ? In fact, closely related numbers
appear in [BoCo] in their description
of the arithmetical symmetry breaking, and in [Jo1], as indices
of subfactors of the hyperfinite factor of type $II_1$.
This suggests interesting questions in the framework
of our program.

\medskip

{\bf 0.6. Elliptic curves as non--commutative spaces.} $\bold{C}$--points
of the elliptic curve $E_{\tau}$ associated with the lattice
$\bold{Z}+\bold{Z}\tau$ can be identified via the exponential map with
$$
\bold{C}/(\bold{Z}+\bold{Z}\tau )\cong \bold{C}^*/(q^{\bold{Z}}),\
q:=e^{2\pi i\tau}.
$$
To treat $E_{\tau}$ as a non--commutative space means to study
the appropriate crossed product of an algebra of functions
on $\bold{C}^*$ with its automorphism group generated
by the shift $z\mapsto qz$. The simplest crossed product
of this type is $A_q^{alg}:=\bold{C}[z,z^{-1}][v,v^{-1}]$
where $vz=qzv,$ but more sophisticated versions
(various completions and their subalgebras) are really
interesting.

\smallskip

There is a lot of results in the theory of lattice models and 
$q$--deformations that can be interpreted in the light  of non--commutative
geometry and function theory of elliptic curves.
It would be worthwhile to review them systematically for two reasons:
first, to get an environment in which elliptic curves
and quantum tori could be treated more or less uniformly,
and second, because the non--commutative setting considerably
enriches even the classical picture.

\smallskip

I will restrict myself by two examples illustrating
these points.

\medskip

{\bf 0.6.1. Semistable bundles and regular modules.} Here we explain
the basic reslut of [BG], generalized in [BEG]. The relevant
crossed product is $A^{form}_q :=\bold{C}((z,z^{-1}))[v,v^{-1}]$
where $\bold{C}((z,z^{-1}))$ denotes the field of formal
Laurent series finite in negative degrees. Consider the following
categories.

\smallskip

{\it Category I.} Its objects are left $A^{form}_q$--modules $M$,
which are finite--dimensional as $\bold{C}((z,z^{-1}))$--spaces
and which satisfy the following {\it regularity} condition:
there exists a free $\bold{C}[[z]]$--submodule $M_0\subset M$
of maximal rank such that $v^{\pm 1}(M_0)\subset M_0.$

\smallskip

Morphisms are usual module homomorphisms. The tensor product
over $\bold{C}((z,z^{-1}))$ extends to a structure of
rigid tensor category.

\smallskip
 
{\it Category II.} Its objects are semistable degree 0 holomorphic
vector bundles over the elliptic curve $E_{\tau}.$
Semistability means that the associated principal bundle
admits a global holomorphic connection,
which is automatically flat. Morphisms and tensor structure
are evident ones.

\smallskip

One of the main results of [BG] consists in the construction
of an equivalence between the categories I and II. This
equivalence is compatible with tensor products.

\smallskip

An interesting comment made in [BG] connects this result with
a problem in the theory of finite difference
equations. 

\smallskip

Consider first differential equations of the form
$$
z\,\frac{d\chi}{dz} = m(z)\,\chi (z)
\eqno(0.9)
$$
where $\chi (z)$ is a column of formal series from
$\bold{C}((z,z^{-1}))$ and $m(z)$ a matrix of such series.
We can try to classify such equations by identifying
those which can be obtained from each other by
a linear transformation of $\chi (z)$. This is equivalent
to the gauge transformation of $m(z)$:
$$
m(z) \mapsto g(z) \,m(z)\, g(z)^{-1}+z\,\frac{dg}{dz}\,g(z)^{-1}.
\eqno(0.10)
$$
It is known that if (0.9) has a regular (Fuchsian) singularity at $z=0$, then
the formal classification coincides with the analytic one,
and the latter is furnished by monodromy around $0$.

\smallskip

Now, a finite difference, or $q$--version of (0.9) is
$$
\chi (qz) = m(z)\,\chi (z)
\eqno(0.11)
$$
and the gauge equivalence (0.10) is replaced by
$$
m(z)\mapsto g(qz) \,m(z)\, g(z)^{-1}.
\eqno(0.12)
$$
We can identify the problem of classification
of the equations (0.11) up to gauge equivalence
with the problem of classification up to isomorphism of $A^{form}_q$--modules
finite--dimensional over $\bold{C}((z,z^{-1})).$
To this end, given $m(z)$, treat it as the matrix
of the operator $v$ in a basis. The Baranovski--Ginzburg
theorem then implies that for regular modules
this classification coincides with
the classification of semistable vector bundles
over $E_{\tau}.$ But the regularity condition
for modules is the standard $q$--version of Fuchsian regularity.
Hence semistable vector bundles over $E_{\tau}$ should be regarded
as a $q$--version of the monodromy data.

\smallskip

It would be important to reconstruct the complete
category of coherent sheaves and/or its derived
category in terms of an appropriate crossed product.
See an interesting discussion in [So], especially 3.3.

\medskip

{\bf 0.6.2. Quantum pentagon identity.} As above, let
$q=e^{2\pi i\tau}$, $\roman{Im}\,\tau >0.$ This time we will consider
the elliptic curve $E_{2\tau}$ represented by an appropriate
completion of the algebra $\bold{C}[u,u^{-1},v,v^{-1}]$
with $uv=q^2vu.$ Put
$$
e_q(t):=\prod_{n\ge 0}(1+q^{2n+1}t).
\eqno(0.13)
$$
If $t$ here is understood as a complex number, we get
one of the standard classical expressions,
for example, occuring in the product formula for the elliptic theta function
$$
\theta_q(t) :=\sum_n q^{n^2}t^n =e_q(t)\,e_q(t^{-1}).
$$
The following noncommutative identities (with $uv=q^2vu$)
are however nonclassical:
$$
e_q(u)\,e_q(v)=e_q(u+v),
\eqno(0.14)
$$
$$
e_q(v)\,e_q(u) =e_q(u)\,e_q(vu)\,e_q(v).
\eqno(0.15)
$$
In view of (0.14), $e_q(t)$ is sometimes called the $q$--exponential
function. 

\smallskip

The second identity (0.15) was proved by Faddeev and Kashaev
in [FK] and called there a quantum version of the
Rogers pentagon identity for the dilogarithm. To explain
this, I remind here the classical version of the Rogers
identity:
$$
L(x)+L(y)-L(xy)=L\left(\frac{x-xy}{1-xy}\right)
+L\left(\frac{y-xy}{1-xy}\right) \, ,
\eqno(0.16)
$$
where
$$
L(x):=L_2(x) +
\frac{1}{2}\,\roman{log}\,(1-x)\,\roman{log}\,x
$$
and
$$
L_2(x):=-\int_0^x \roman{log}\,(1-z)\,\frac{dz}{z}=\sum_{n\ge 1}
\frac{x^n}{n^2}.
$$
As $q$ tends to the cusp 1, we have a classical asymptotic expansion
in $\tau$ for the logarithm of $e_q(t)$ for which we write a few first terms
in the exponentiated form
$$
e_q (t)=\frac{1}{\sqrt{1+qt}} \, \roman{exp}\,(L_2(-t)/4\pi i\tau )\,
(1+ O(\tau )).
\eqno(0.17)
$$
It remains to combine (0.15) and (0.17). This is not quite
straightforward. Faddeev and Kashaev argue that
an appropriate infinite--dimensional representation of the
commutation relations $uv=q^2vu$ and the corresponding notion
of the symbol of an operator in this representation
produce (0.16).

\smallskip

It is remarkable and promising that (0.15) looks much neater than
(0.16) and shows that (0.16) is a boundary reflection
of a phenomenon which is both more global
and essentially noncommutative.


\bigskip

\centerline{\bf \S 1. Pseudolattices, quantum tori,}

\smallskip

\centerline{\bf and Real Multiplication}

\medskip

{\bf 1.1. Category of pseudolattices $\Cal{P}\Cal{L}$.} By definition, {\it a pseudolattice}
(of rank 2) is a quadruple $(L,V,j,s)$, where 
$L$ is a free abelian group of rank two, $V$ is
an one--dimensional complex space, $j:\,L\to V$ is an injective homomorphism
whose image lies on a real line, 
and finally $s$ is an orientation of this line. Since this line
contains $0$, a choice of $s$ defines the notion of positive
and negative halves of it. Clearly, this line is the topological
closure of $j(L).$

\smallskip

 {\it A strict morphism}
of pseudolattices $(L^{\prime},V^{\prime},j^{\prime},s^{\prime})\to (L,V,j,s)$ is a commutative diagram
$$
\CD 
L^{\prime} @>j^{\prime}>> V^{\prime}\\
@V\varphi VV  @VV\psi V  \\
L @>>j> V  \\
\endCD
\eqno(1.1)
$$
in which $\varphi$ is a group homomorphism, and
$\psi$ is a $\bold{C}$--linear map, which transforms
the orientation $s^{\prime}$ to $s$. Clearly, 
$\phi$ and $\psi$ uniquely determine each other. Moreover, such a strict morphism
is a strict isomorphism iff both $\phi$ and $\psi$ are isomorphisms.

\smallskip

Omitting the condition that $\psi$ respects orientations,
we get the notion of {\it weak morphism.}

\smallskip

As with lattices, several simple observations will help us to
clarify the structure of this category.

\medskip

(i) The orientation $s$ makes  $L$ a totally ordered group:
by definition, $l>m$ iff $j(l-m)$ lies in the $s$--positive half--line.
Choosing a basis $(l_1,l_2)$ in $L$ and 
taking $j(l_2)$ as the base vector of $V$, we see that
in any strict isomorphism class of pseudolattices 
one can find a representative given by $j:\,\bold{Z}^2\to\bold{C}$
such that $j(0,1)=1$, $j(1,0):=\theta$ is an irrational real number.
The remaining piece of data is
the sign $\varepsilon =\pm 1$ such that
$l>0$ iff $\varepsilon j(l)>0.$

\smallskip

Let us denote this pseudolattice $(L_{\theta}, \varepsilon ).$
Then any non--zero strict morphism $(L_{\theta^{\prime}}, \varepsilon^{\prime} )\to (L_{\theta}, \varepsilon )$ is represented
by a non--degenerate matrix 
$$
g=\left(\matrix a&b\\c&d\endmatrix\right)\in M\,(2,\bold{Z})
$$
such that
$$
 \theta^{\prime} =\frac{a\theta +b}{c\theta +d},\quad
\roman{sgn}\,(c\theta +d)=\varepsilon\varepsilon^{\prime}.
\eqno(1.2)
$$ 
This $g$ is obtained by writing  $\varphi$ in (1.1)
as the right multiplication of a row by a matrix. 

\smallskip

If $g\in GL(2,\bold{Z})$, (1.2) is an isomorphism. Since we can 
replace $g$ by $-g$ without violating the first condition in (1.2),
two lattices are strictly isomorphic iff they are weakly isomorphic.

\medskip 

(ii) We can choose a positive basis in $L$. This shows that
any pseudolattice is isomorphic to some $(L_{\theta}, \varepsilon =1)$
with irrational real $\theta$  which can be even taken in (0,1). 
We will denote it simply $L_{\theta}.$
Two such pseudolattices are isomorphic iff their invariants
$\theta$ lie in the same $PGL(2,\bold{Z})$--orbit,
that is, their continued fraction expansions coincide starting 
from some place. Thus set--theoretically, the moduli space
of the isomorphism classes of pseudolattices is
$$
(PGL(2,\bold{Z})\setminus \bold{P}^1(\bold{R}))\,\setminus \, \{\roman{cusp}\}
\eqno(1.3)
$$
where the cusp is the orbit of rational numbers.

\medskip

(iii) Weak endomorphisms of a pseudolattice $L$
(we omit other structures in notation if there is no
danger of confusion) form a ring
$w$-$\roman{End}\,L$ (w stands for weak),
with componentwise addition of $(\phi, \psi)$ and
composition as multiplication.  It contains $\bold{Z}$
and comes together with its embedding in $\bold{R}$
as $\{ a\in\bold{R}\,|\, a j(L)\subset j(L) \}.$
The non--negative part of this ring is the semiring
$\roman{End}\,L$.

\proclaim{\quad 1.1.1. Lemma} (a) w-$\roman{End}\,L\ne \bold{Z}$
iff there exists a real quadratic subfield $K$ of $\bold{R}$ 
such that $L$ is isomorphic to a pseudolattice contained in $K.$

\smallskip

(b) If this is the case, we will say that
$L$ is an RM pseudolattice. Denote by $O_K$ the ring of integers of
$K$. There exists a unique integer $f\ge 1$ (conductor) such that
w-$\roman{End}\,L =\bold{Z}+fO_K=:R_f,$ and $L$ is
a projective module of rank 1 over $R_f$. 
\smallskip
The module $L$ is endowed with a total ordering
respected by $\roman{End}\,L$.

\smallskip

Every $K$, $f$
and a ordered projective module over $R_f$ come from a lattice.

\smallskip

(c) If pseudolattices $L$ and $L_1$ have the same
$K$ and $f$, they are isomorphic if and only if
their classes in the Picard group  $\roman{Pic}\,R_f$ coincide. 
\endproclaim

\medskip

Unlike the case of lattices, the automorphism group
of a pseudolattice is always infinite, it is isomorphic
to $\bold{Z}\times \bold{Z}_2.$ 

\smallskip

For RM pseudolattices embedded in one and the same real
quadratic field $K$, we will say that an
isomorphism $L\to L_1: l\mapsto al_1,\,a\in K$,
is an isomorphism {\it in a narrow sense}, if $N_{k/\bold{Q}}(a)>0.$

\medskip

{\bf 1.2. Two--dimensional quantum tori.} We now want to define
analogs of elliptic curves for pseudolattices, that is, some geometric
objects representing quotients $V/j(L)$
where $(L,V,j,s)$ is a pseudolattice.

\smallskip

Choosing $L_{\theta}$ as a representative of the
respective isomorphism class, we can naively replace
$\bold{C}/(\bold{Z}+\bold{Z}\theta )$ by
$\bold{C}^*/(e^{2\pi i\theta})$ (``Jacobi uniformization''),
and then interpret the last quotient
as an ``irrational rotation algebra'', or
two--dimensional quantum torus $T_{\theta}.$ We recall
that this torus
is (represented by) the universal $C^*$--algebra $A_{\theta}$
generated by two unitaries $U,V$ with the commutation rule
$UV=e^{2\pi i\theta}VU.$ A choice of such generating
unitaries is called {\it a frame}; it is not unique.

\smallskip

The next task is to define morphisms between these quantum tori,
with properties that would allow us to imitate the
framework of 0.3. Already isomorphisms present a problem:
we want fractional linear transforms (1.2) to produce
isomorphic quantum tori. M.~Rieffel's seminal
discovery was that to this end we should consider
Morita equivalences between appropriate categories of modules 
as isomorphisms between the tori themselves.
Morita equivalences are essentially given by tensor multiplication
by a bimodule.
Taking this lead, we will formally introduce the general
Morita morphisms of associative rings, stressing those
traits of the formalism that play a central role in the
structure theory of quantum tori (of arbitrary dimension).
Our presentation also  prepares ground for
introducing versions of quantum tori with more 
algebraic geometric flavor.

\medskip

{\bf 1.3. Morita category.}  Let $A,B$ be two associative rings.
A {\it Morita morphism} $A\to B$
by definition, is the isomorphism class of a bimodule 
${}_AM_B$, which is projective and finitely generated
separately as module over $A$ and  $B$.

\smallskip

The composition of morphisms is given by the tensor product
${}_AM_B\otimes {}_BM^{\prime}_C,$ or ${}_AM\otimes {}_BM^{\prime}_C$
for short.

\smallskip  

If we associate to ${}_AM_B$
the functor
$$
\roman{Mod}_A\to \roman{Mod}_B:\ 
N_A\mapsto N\otimes_AM_B,
$$
the composition of functors will be given by the tensor product,
and isomorphisms of functors will correspond to the isomorphisms
of bimodules.

\smallskip

We imagine an object $A$ of the (opposite) Morita category  as a  
noncommutative space, right $A$--modules as sheaves on this
space, and the tensor multiplication by ${}_AM_B$
as the pull--back functor, in the spirit of
A.~Rosenberg's program [Ro2].
We have chosen to work with right modules, but passing to 
the opposite rings allows one to reverse
left and right in all our statements.

\smallskip

Two bimodules ${}_AM_B$ and ${}_BN_A$ supplied with two
bimodule isomorphisms ${}_AM\otimes_BN_A \to {}_AA_A$ 
and ${}_BN\otimes_AM_B\to {}_BB_B$ define
mutually inverse Morita isomorphisms (equivalences) between
$A$ and $B$. The basic example
of this kind is furnished by $B=\roman{Mat}\,(n,A)$,
$M={}_AA^n{}_B$ and $N={}_BA^n{}_A$.

\smallskip

We will now briefly summarize Morita's theory.

\medskip

(A) {\it Characterization of functors $S:\,\roman{Mod}_A\to \roman{Mod}_B$ 
of the form $N_A\mapsto N\otimes_AM_B.$} They are precisely functors
satisfying any of the two equivalent conditions:

\smallskip

(i) {\it $S$ is right exact and preserves direct sums.}

\smallskip

(ii) {\it $S$ admits a right adjoint functor 
$T:\,\roman{Mod}_B\to \roman{Mod}_A$ (which is then
naturally isomorphic to $\roman{Hom}_B(M_B,*)$).}

\smallskip

We will call such functors {\it continuous}.

\medskip

(B) {\it Characterization of continuous functors $S$
such that $T$ is also continuous and $ST\cong 1.$} 
Let $S$ be given by ${}_AM_B$ and $T$ by ${}_BN_A$. 
Then $M\otimes_BN\cong {}_AA_A.$ Moreover, in this case

\smallskip

(iii) {\it $M_B$ and ${}_BN$ are projective.}

\smallskip

(iv) {\it ${}_AM$ and $N_A$ are generators.}

\smallskip

In particular, equivalences $\roman{Mod}_A\to \roman{Mod}_B$ are 
automatically continuous. Hence any pair of
mutually quasi--inverse equivalences must be given
by a couple of biprojective bigenerators as above.

\smallskip

(C) {\it Finite generation and balance.}
Any right module $M_B$ can be considered as a bimodule
${}_AM_B$ where $A=B':=\roman{End}_B(M_B).$ We can then
similarly produce the ring $B''=A':= \roman{End}_A({}_AM).$
Module $M_B$ is called {\it balanced} if $B''=B.$
Similarly, one can start with a left module. With this
notation, we have:

\smallskip

{\it (v) $M_B$ is a generator iff ${}_{B'}M$ is balanced
and finitely generated projective.}

\medskip

Properties (i)--(v) can serve as a motivation
for our definition of the Morita category above.

\medskip

{\bf 1.4. Projective modules, idempotents, traces, and $K_0$.}
Projective right $A$--modules up to isomorphism
are exactly ranges of idempotents in various
matrix rings $\roman{Mat}\,(n,A)$ acting from the left upon
(column) vector modules $A^n$. Morphisms between
such modules are also conveniently described
in terms of these idempotents. The following (well known)
Proposition summarizes the relevant information
in the form convenient for us.

\smallskip

We prefer to work with all $n$ simultaneously.
So we will denote by $\Cal{M}A$ the ring of infinite matrices
$(a_{ij})\, i,j\ge 1,\,a_{ij}\in A,$ $a_{ij}=0$ for $i+j$
big enough (depending on the matrix in question). Notice that
$\Cal{M}A$ is not unital even if $A$ is. Similarly, denote
by $A^{\infty}$ the left $\Cal{M}A$--module of infinite 
columns $(a_i),\, i\ge 1,$ with coordinates in $A$ and such that
$a_i=0$ for large $i.$

\smallskip

Denote by $Pr_A$ the category of finitely generated right $A$--modules.
Denote by $pr_A$ the category whose objects are
projectors (idempotents) $p\in \Cal{M}A,\,p^2=p$, whereas morphisms are
defined by
$$
\roman{Hom}\,(p,q):=q\,\Cal{M}A\,p
$$
with the composition induced by the multiplication
in $\Cal{M}A$.

\smallskip

There is a natural functor $pr_A\to Pr_A$ defined on objects
by $p\mapsto pA^{\infty}$. In order to define it on
morphisms, we remark that morphisms $pA^{\infty}\to qA^{\infty}$
can be naturally described by matrices in the following way.
Clearly, $p A^{\infty}$ contains the columns
$p_k$ of $p$ which generate $p A^{\infty}$ as a right $A$--module.
We can then apply any $\varphi :\, p A^{\infty}\to q A^{\infty}$
to all $p_k$ and arrange the resulting vectors into a matrix
$\Phi\in \Cal{M}A$ with $k$--th column $\varphi (p_k).$
One checks that $\Phi p=\Phi$, and since also
$q \Phi =\Phi$, we have $\Phi\in q\Cal{M}Ap$. Conversely,
any such matrix determines a unique morphism $pA^{\infty}\to qA^{\infty}$.

\medskip

\proclaim{\quad 1.4.1. Proposition} (a) The functor
$pr_A\to Pr_A$ described above is an equivalence of categories.

\smallskip

(b) $p A^{\infty}$ is isomorphic to $q A^{\infty}$ iff
there exist $X,Y\in \Cal{M}A$ such that $p=XY, q=YX$
(von Neumann's equivalence of idempotents).
Hence in this case $p-q\in [\Cal{M}A,\Cal{M}A].$
\endproclaim 

We have already checked most of the statements implicit in (a).
As for (b), consider two mutually inverse isomorphisms 
$p A^{\infty} \to q A^{\infty}$ and $q A^{\infty} \to p A^{\infty}$.
Assume that the first one sends columns of $p$ to the columns
of $qBp$ whereas the second one sends columns
of $q$ to the columns of $pCq$. Writing that their compositions
send $p$ to $p$ and $q$ to $q$, we see that one can take $X=pCq$,
$Y= qBp.$ Conversely, if $p=XY, q=YX$, then also 
$p=(pXq)(qYp)$, $q=(qYp)(pXp)$, so that the matrices $pXq$ and
$qYp$ determine mutually inverse isomorphisms of $p A^{\infty}$   
and $q A^{\infty}$.
This completes the proof.

\medskip

It is also convenient to introduce the parallel formalism
for left projective $A$--modules: here we consider
the right $\Cal{M}A$--module $A_{\infty}$ of 
rows $(a_i),\, i\ge 1,\, a_i=0$ for large $i$, and map a projector $p$ to the left $A$--module
$A_{\infty}p.$

\smallskip

Replacing $A$ by the opposite ring $A^{op}$ switches these two
constructions. 

\smallskip

{\it A trace} of $A$ is any homomorphism of additive groups $t:\,A\to G$
vanishing on commutators; by definition, it
factors through the universal trace $A\to A/[A,A]$.
Combining it with the matrix trace, we get
its canonical extension to $\Cal{M}A$. From the Proposition 1.4.1 (b)
it follows that $t(p)$ depends only on the
isomorphism class of $pA^{\infty}.$
The class $p\,\roman{mod}\, [A,A]$
 is called {\it the Hattori--Stallings rank} of $p A^{\infty}$.

\smallskip

We define $K_0(A)$ as the Grothendieck group
of $Pr_A$. If $N_A\in Pr_A$,
$[N_A]$ denotes its class in $K_0(A).$
If $N_A$ is the range of an idempotent $p$
and $t$ is a trace, $t(p)$ depends only on $[N_A]$ and is
additive on exact triples, hence 
$t$ becomes a homomorphism of $K_0(A)$
(it is called dimension in the theory of von Neumann algebras).

\smallskip

 The crucial role
of traces
in the theory of quantum tori (and in more general
functional analytic situations) is explained
by the fact that for irrational tori
projective modules are exactly classified by the
value of the (unique) normalized trace of the
respective projector (Rieffel).

\smallskip

The following simple Lemma on traces is a useful technical tool.
We assume in it that we work with algebras
over a ground field.

\smallskip

\proclaim{\quad 1.4.2. Lemma} Consider two unital algebras $A,B$ 
and an $A$--$B$--bimodule ${}_AM_B$ which is projective as a module
over $A$ and over $B$. 
Assume that the (dual) space of traces $A/[A,A]$ of $A$ is one--dimensional,
whereas that of $B$ is $\ge 1$--dimensional, and that $1\notin [A,A]$.
Choose non--zero traces $t_A$ and $t_B$. Then there
exists such a constant $c$ that for any $N_A\in Pr_A$
we have
$$
t_B([N\otimes_AM_B])=c\,t_A([N_A])).
\eqno(1.4)
$$
The value of this constant is obtained by putting $N_A=A_A$ in
(1.4):
$$
c=t_B([{}_AM_B])t_A([A_A])^{-1}.
\eqno(1.5)
$$
\endproclaim
\smallskip

{\bf Proof.} Let the  modules  $N_A$, $M_B$
be given as ranges of idempotents $q\in \Cal{M}A$, $p\in \Cal{M}B$
respectively.

\smallskip

Put $A_1:=\roman{End}_B(M_B)$ and identify 
$p\Cal{M}Bp$ with $A_1$ as above. The structure of left 
$A$--module on $M_B$ is given
by a ring homomorphism $\varphi:\,A\to A_1.$
The trace $t_B$ induces a trace on $p\Cal{M}Bp$
and thus a trace $t_1$ on $A_1.$ In turn, $t_1$
induces via $\varphi$ a trace on $A$, and since
the latter is unique, there exists a constant $c\ne 0$
such that we have 
identically $t_1(\varphi (a))=c\,t_A(a)$ for all $a\in A.$
By the associativity of tensor multiplication, we have
$$
N\otimes_AM_B=(N\otimes_AA_1)\otimes_{A_1}M_B.
$$
As a right $A_1$--module, $N\otimes_AA_1$ is isomorphic to $q_1A_1^{\infty}$ where
$q_1=\varphi (q).$
We have $q_1^2=q_1$ and since $q_1\in p\Cal{M}Bp$,
$q_1p=pq_1=q_1$.  Hence finally $N\otimes_AM_B$ as a right $B$--module
is isomorphic to $q_1B^{\infty}$. Thus
$$
t_B([N\otimes_AM_B])=t_B(q_1)=t_{A_1}(q_1)=c\,t_A(q)=
c\,t_A([N_A]).
$$

\medskip

{\bf 1.5. Involutions and scalar products.} Assume now that
$A$ is endowed with an additive
(linear or antilinear in the case of algebras) involution $a\mapsto a^*$, $(ab)^*=b^*a^*,$
$a^{**}=a.$
It extends to matrix algebras: $(B^*)_{ij}:=B_{ji}^*.$
Similarly, it extends to $A^{\infty}\to A_{\infty}$
and $ A_{\infty}\to A^{\infty}$, compatibly
with the module structures.

\smallskip

In such a context, it makes sense to consider only those
projective modules which are ranges of {\it projections},
that is, $*$--invariant idempotents. In fact, in the case
of $C^*$--algebras the resulting subcategory of projective modules
is equivalent to the full category, because
{\it every idempotent is von Neumann equivalent to
a projection.} In fact, if $p$ is an idempotent,
then 
$$
P:=pp^*\,[1-(p-p^*)^2]^{-1}
$$
is an equivalent projection. The $C^*$--structure is used
to ensure the invertibility of $1-(p-p^*)^2$; the rest is pure algebra:
see e.~g. [Da], IV.1.

\smallskip

Taking into account the involution, we  get
an additional structure on our modules and bimodules
consisting of scalar products and identities relating them.
This is a simple but important formalism made explicit by M.~Rieffel.

\medskip

\proclaim{\quad 1.5.1. Lemma} Let $M_B$ be a projective
module (over a ring with
involution $B$) isomorphic to $p\Cal{M}B$ with $p^*=p.$ Put $A=\roman{End}_B(M_B)$,
identify this ring with $p\Cal{M}Bp$ as above, and consider $M$ as an $A$--$B$ bimodule. The involution on $pBp$ is induced by that on $B$.

\smallskip

Define two scalar products 
${}_A\langle*,*\rangle:\,M\times M\to A$ and 
$\langle*,*\rangle_B:\,M\times M\to B$:
$$
  {}_A\langle pb,pc\rangle :=(pb)(pc)^*=pbc^*p\in p\Cal{M}Bp=A,
\eqno(1.6)
$$
$$
  \langle pb,pc\rangle_B :=(pb)^*pc=b^*pc\in B .
\eqno(1.7)
$$
Then the following identities hold, in which $l,m,n\in M$, $a\in A,b\in B$:
$$
{}_A\langle m,n\rangle^* ={}_A\langle n,m\rangle ,\quad
\langle m,n\rangle_B^* =\langle n,m\rangle_B ,
\eqno(1.8)
$$
$$
{}_A\langle am,n\rangle =a {}_A\langle m,n\rangle,\quad
{}_A\langle m,an\rangle ={}_A\langle m,n\rangle\, a^*,
\eqno(1.9)
$$
$$
\langle mb,n\rangle_B=b^*\langle m,n\rangle_B,\quad
\langle m,nb\rangle_B =\langle m,n\rangle_Bb ,
\eqno(1.10)
$$
$$
{}_A\langle l,m\rangle n= l \langle m,n\rangle_B.
\eqno(1.11)
$$
\endproclaim

\smallskip

We omit the checks which are straightforward.

\medskip

{\bf 1.6. Rieffel's projections.} As we will see in \S 3,
over toric $C^*$--algebras
many bimodules ${}_AM_B$  are constructed directly,
by inducing them from a Heisenberg representation,
and the scalar products with the properties summarized
in the Lemma 1.4.2 are introduced by an ad hoc
formula.

\smallskip

In this case it is useful to know that, conversely,
projections can be produced from such a setup.
The following Lemma due to Rieffel ([Ri3]) furnishes them.

\medskip

\proclaim{\quad 1.6.1. Lemma} Assume that ${}_AM_B$
is a bimodule over two rings with involution, endowed with
two scalar products satisfying the formalism (1.8)--(1.11).
Let $m\in {}_AM_B.$

\smallskip

(a) If $m\langle m,m\rangle_B=m$, then $p:={}_A\langle m,m\rangle$
is a projection in $A$.

\smallskip

(b) Conversely, assume that from ${}_A\langle n,n\rangle =0$
it follows that $n=0$. In this case, if $p$ as above is a projection,
then $m\langle m,m\rangle_B=m$.
\endproclaim 

\smallskip

{\bf Proof.} (a) Using (1.6) and (1.8), we obtain
$$
p^2={}_A\langle m,m\rangle{}_A\langle m,m\rangle =
{}_A\langle {}_A\langle m,m\rangle m,m\rangle =
{}_A\langle m\langle m,m\rangle_B ,m\rangle
={}_A\langle m,m\rangle =p.
$$
From (1.3) it follows that $p^*=p.$

\smallskip

(b) Conversely, if $p$ is a projection, then we get similarly
$$
{}_A\langle m\langle m,m\rangle_B -m , m\langle m,m\rangle_B-m\rangle =0.
$$
\smallskip

This completes the proof.

\medskip

Rieffel also remarks that if $x\in {}_AM_B$ is such an
element that one can construct an invertible
*--invariant square root $\langle x,x\rangle_B^{1/2}$,
then $m:=x\,\langle x,x\rangle_B^{-1/2}$ satisfies 1.6.1(a).

\smallskip

F.~Boca in [Bo2] takes for $x$ a Gaussian element 
in the relevant Heisenberg module. Then $\langle x,x\rangle_B$
turns out to be a quantum theta in the sense of [Ma3]. We develop
this remark in \S 3 for multidimensional case. 

\smallskip

The net result is that  we have a supply of explicit projections
in toric algebras which in the notation of \S 3 below are given by the formulas
$$
p_T:={}_D\langle f_T\,\Theta_{D^!}^{-1/2}, f_T\,\Theta_{D^!}^{-1/2}\rangle
\in C(D,\alpha ) .
$$
For notation, see  (3.20), Theorem 3.7, and section 3.3.

\smallskip

This formula (and its generalizations) relates the representation theory of
quantum tori to the theory of quantum thetas which has a distinct flavor
of non--commutative {\it algebraic} geometry.
Notice however that the existence of $\Theta_{D^!}^{-1/2}$
is not established in full generality.

\smallskip

The whole formalism sketched in 1.3--1.6.1 is a simple
algebraic version of some basic machinery in the theory
of von Neumann and $C^*$--algebras. In particular,
see [Co1] and [Jo2] and original papers by A.~Connes and A.~Wassermann
who overcame some highly nontrivial complications 
arising in the operator context.

\medskip

{\bf 1.6.2. Morita category and two--dimensional quantum tori.}
By definition, two--dimensional quantum tori are objects
of the category $\Cal{Q}\Cal{T}$ whose morphisms are
isomorphism classes of projective bimodules ${}_AM_B$
corresponding to  projections, so that
the formalism of the previous subsections is readily
applicable. In particular, $A_{\theta}$
has a unique trace $t_A$ which is normalized
by the condition $t_A(1)=1$ and which vanishes
on any frame.

\smallskip

We will see in 1.7 below that there is a functorial
correspondence between $\Cal{Q}\Cal{T}$ and
pseudolattices which is fairly similar to the
correspondence between elliptic curves and lattices.
In particular, Real Multiplication of pseudolattices
is reflected in $\Cal{Q}\Cal{T}$.

\smallskip

In order to achieve arithmetical
applications of Real Multiplication, one 
has to find still smaller rings and modules,
perhaps finitely generated in an algebraic sense
and admitting models over rings of algebraic integers.
Their definition remains the central unsolved problem in our approach.
Since the points of finite order $m$ on an elliptic curve $E/K$
are in fact points of a finite group scheme over $K$ acting
upon $E$, it is conceivable that in the
$C^*$--world the relevant
finite objects should be seeked among weak Hopf algebras
(or weak quantum groupoids)
acting upon $C^*$--algebras: see recent reports
[NiVa], [KaNi1], [KaNi2], and the references quoted therein.

\smallskip

The famous paper [Jo1] shows how a spectrum of algebraic
numbers can be generated from such a setting.
Jones's discrete spectrum
of indices of subfactors is $\{4\,\roman{cos}^2\,\dfrac{\pi}{n}\,|\,n\ge 3\}$,
whereas Stark's numbers in the cyclotomic case are
$4\,\roman{sin}^2\,\dfrac{\pi m}{n}$. Both generate the maximal
real subextension of $\bold{Q}^{ab}.$ 

\smallskip

Is this only a coincidence?

\medskip

Returning to the $C^*$ (or smooth) context,
notice in conclusion that a bimodule ${}_AM_B$ can be treated
as an $A\otimes B^{op}$--left module (completed
tensor product). If it were projective, we could classify
bimodules for toric $A,B$ using the fact that $A\otimes B^{op}$
is again toric: their invariants would come from
$K_0(A)\otimes K_0(B^{op})$ (the ``trivial part'')
and from $K_1(A)\otimes K_1(B^{op})$ (the really interesting
correspondences). However, intuitively it seems clear
that such bimodules are much smaller than projective
modules because they are separately $A$-- and $B^{op}$--projective
and hence, like Rieffel's elementary modules, should be
realisable in functions of $\roman{dim}\,A=\roman{dim}\,B$
variables, whereas $A\otimes B^{op}$--projective modules
are realizable only by functions
of the doubled number of variables.

\smallskip

Therefore several questions arise about a possible
extension of the classification theory of modules.

\smallskip

{\it Question.} Is any Morita morphism bimodule  a maximal
quotient of a unique projective  $A\otimes B^{op}$--module ?

\smallskip

More generally, toric projective modules can have nontrivial
maximal quotients, like in the situation with highest weight
and Verma modules. One should seek for canonical
projective resolutions of such modules.

\smallskip

The algebraic machinery might be connected with the fact that
$A\otimes B^{op}$ contains large commutative subalgebras,
so that a module can be decomposed
according to their characters. E.~g. if $A,B$ are two--dimensional
quantum tori, $A\otimes B^{op}$ contains two--dimensional
classical tori, and prescribing their characters
may produce the interesting quotients.

\smallskip

{\it Question.} Can one find a description of the
derived category of perfect complexes over toric algebras?

\medskip
{\bf 1.7. Two functors relating $\Cal{Q}\Cal{T}$ to $\Cal{P}\Cal{L}$.}
We start with defining a functor $K:\,\Cal{Q}\Cal{T} \to \Cal{P}\Cal{L}.$
Let the torus $T$ be represented by an algebra $A$. On objects, we put:
$$
K(T)=(L_A,V_A,j_A,s_A).
\eqno(1.12)
$$
Here
$L_A:=K_0(A)$, the $K_0$--group
of the category of right projective $A$--modules
(as above, given by projections in finite matrix algebras over $A$); 
$V_A$ is the target group of the universal trace on $A$,
that is, the quotient space of $A$ modulo the completed
commutator subspace $[A,A]$. Furthermore, 
$j_A=t_A:\,K_0(A)
\to V_A$
is this universal trace extended to matrix algebras; its value
on the class of a module, as we already explained, is its value
at the respective projection.
Finally $s_A$ is taken in such a way
that positive elements in $K_0(A)$ become represented
by the classes of actual (not virtual) projective modules.

\smallskip

On morphisms, we define directly
the left vertical arrow of the respective diagram (1.1):
$$
K({}_AM_B)([N_A]):=[N\otimes_AM_B].
\eqno(1.13)
$$
The existence of the right vertical arrow follows
from the Lemma 1.4.1.

\medskip

\proclaim{\quad 1.7.1. Theorem} (a) The family of maps (1.12), (1.13) can be uniquely completed to a functor $K:\, \Cal{Q}\Cal{T}\to \Cal{P}\Cal{L}$.

\smallskip

(b) This functor is essentially surjective on objects and (strict) morphisms.

\smallskip

(c) Assume that two bimodules ${}_AM_B$ and ${}_AM_B^{\prime}$
considered as morphisms in $\Cal{Q}\Cal{T}$ become equal after
applying $K$. Put $A_1:=\roman{End}_B(M_B)$ and consider
${}_AM_B$ as an $A_1$--$B$ bimodule ${}_{A_1}M_B$ .

\smallskip
There exist
two ring homomorphisms $\varphi, \psi :\,A\to A_1$ such that
if one considers $A_1$ as an $A$--$A_1$ bimodule ${}_{\varphi}A_{1,A_1}$,
(resp.   ${}_{\psi}A_{1,A_1}$) using $\varphi$ (resp. $\psi$)
to define the left action, 
and the ring structure of $A_1$ to define the right action,
one obtains
$$
{}_{\psi}A_1\otimes_{A_1}M_B^{\prime}\cong {}_{\varphi}A_1\otimes_{A_1}M_B
\eqno(1.14)
$$
as $A$--$B$--bimodules.

\smallskip

In particular, if $\otimes {}_AM_B$
and $\otimes {}_AM_B^{\prime}$ produce Morita equivalences,
these functors differ by an automorphism
of the category $\roman{Mod}_A$ which is induced by
an automorphism of the ring $A$.

\endproclaim

{\bf Comments.} This result should be compared to
the easy Theorem 0.3.1 which provides the geometric
basis of the Complex Multiplication. 
The statement about quantum tori sounds less neat,
however in 1.7.2 we will complement it by the construction
of a functor in the reverse direction
defined only on isomorphisms,
which should suffice for the envisioned applications
to Real Multiplication.

\smallskip

{\bf Proof.} (a) Lemma 1.4.2 shows that, after passing
to traces, (1.13) becomes the multiplication
by a positive number representing
a (strict) morphism of pseudolattices $K(m).$ 
The compatibility with the composition of morphisms
is straightforward.

\smallskip

(b) It remains to
establish the following three facts.

\smallskip

(i) {\it Every object of $\Cal{P}\Cal{L}$ is isomorphic
to an object lying in the image of $K$.}

\smallskip

In fact, the pseudolattice denoted $(L_{\theta},1)$ in 1.1(b) 
is isomorphic to $K(A_{\theta})$ where $A_{\theta}$
is the respective rotation algebra. This is the main result
of the theory, due to Connes, Rieffel, Pimsner--Voiculescu, Elliott.
It is worth recalling here one of the several known strategies for proving it (cf. [Da], Ch. VI). 

\smallskip

First, one checks that
for any $\alpha\in [0,1]\,\cap\,\bold{Z}+\bold{Z} \theta$
there exists a projection $p_{\alpha}\in A_{\theta}$
with the normalized trace $\tau (p_{\alpha})=\alpha$. Using functional calculus, one can
directly construct such projections of the form
$f(U)V+g(U)+h(U)V^*$ (Rieffel--Powers, see [Da], p. 171.)
It follows than $\tau (K_0(A_{\theta}))\supset \bold{Z}+\bold{Z} \theta .$

\smallskip

Second, one shows that $A_{\theta}$ can be embedded
into an approximately finite algebra $\Cal{A}_{\theta}$ which is
the completed inductive limit of $A_{p_n/q_n}$, where 
$p_n/q_n$ are consecutive convergents to $\theta$.
This embedding allows one to calculate $\tau (K_0(A_{\theta}))$
as the inductive limit of ordered groups $\sigma (K_0(\Cal{A}_{p_n/q_n}))$,
and this inductive limit is explicitly identified with 
$\bold{Z}+\bold{Z} \theta .$

\smallskip

This last argument can be read as a weak continuity
property of $A_{\theta}$ with respect to $\theta$
varying in the set of cusps. In 1.8 below, we will
discuss in what sense $A_{\theta}$ can be regarded as
a limit of $E_{\tau}$ when $\tau$ tends to $\theta$
from the upper half--plane.

\smallskip

(ii) {\it Every morphism $K (T_{\theta})\to K(T_{\theta^{\prime}})$ in  $\Cal{P}\Cal{L}$ is of the form $K(m)$ where
$m$ is the tensor multiplication by an appropriate bimodule.}

\smallskip

Clearly, it suffices to choose a generating family of morphisms in
$\Cal{P}\Cal{L}$ (such that any morphism is a composition
of members of this family) and to show that each generator
can be lifted to $\Cal{Q}\Cal{T}.$

\smallskip

Any morphism of pseudolattices restricted
upon $L$--components is a composition
of an injection and an isomorphism (respecting ordering);
moreover, this restriction uniquely determines it.
Any injection can be decomposed into product
of two injections with cyclic quotients.

\smallskip

Isomorphisms between pseudolattices $L_{\theta^{\prime}}\to L_{\theta}$
can be decomposed into
a sequence of transformations of the form 
$\theta\mapsto -\theta$, $\theta\mapsto \theta +1$,
$\theta\mapsto \theta^{-1}.$ The map $(U,V)\mapsto (V^{\prime},U^{\prime})$
produces an isomorphism $T_{\theta}\to T_{-\theta}$,
whereas $T_{\theta}$ and $T_{\theta+n}$ are obviously
the same. The only non--trivial problem is to find
a Morita equivalence $T_{\theta}\to
T_{\theta^{-1}}$. Its solution was given in [Co2] and generalized to
multidimensional tori in a series of works of Rieffel and his
collaborators, see [Ri5], [RiSch].

\smallskip

Alternatively, in [CoDSch] one can find a direct description of 
a bimodule furnishing a Morita equivalence between
$T_{\theta}$ and $T_{\theta^{\prime}}$, where $\theta$ and $\theta^{\prime}$
are related by a transformation from $PGL(2,\bold{Z}).$
We will reproduce it in 1.7.2 below.

\smallskip

It remains to treat the case of embedding of pseudolattices. 
Now choose $n>0$ and consider the embedding of toric algebras
$B:=A_{n\theta}\hookrightarrow A:=A_{\theta}$ where in self--evident
notation $U_B=U_A^n, V_B=V_A.$ For ${}_AM_B$ take the
bimodule  ${}_AA_B.$ It is free of rank
1 (resp. $n$) as $A$--  (resp. $B$--) module.
Denote by $t_A, t_B$ the normalized traces
(taking value $1$ on $1$). Then the constant
(1.5) is $n$, so that the tensor multiplication
by ${}_AM_B$ produces the morphism of pseudolattices
$L_{\theta} \mapsto L_{n\theta}$: $\theta\mapsto n\theta$, $1\mapsto n.$   
Clearly, any embedding of pseudolattices with cyclic
quotient is isomorphic to such one.

\smallskip

(iii) {\it If $K(m)=K(m^{\prime})$, the respective
$B$--modules $M_B$, $M_B^{\prime}$ are isomorphic.}

\smallskip

In fact, from (1.4) and (1.5) it follows that the $B$--traces
of them coincide, and for two--dimensional irrational tori
this means that they are isomorphic.

\smallskip

The remaining argument is straightforward. 
Choose and fix an isomorphism of $M$ and $M^{\prime}$
as $B$--modules. Actions of $A$ upon $M$ and $M^{\prime}$
correspond to two different homomorphisms
$A\to A_1$. This is the essence of (1.14).

\smallskip

This finishes the proof. 

\medskip

{\bf 1.7.2. The functor $E:\, \Cal{P}\Cal{L}_{iso}\to \Cal{Q}\Cal{T}_{iso} .$}
In this section we rephrase the content of \S 2 of the recent
preprint [DiSch].

\smallskip

Denote by $\Cal{P}\Cal{L}_{iso}$ the category whose objects are
pseudolattices $L_{\theta}=\bold{Z}+\bold{Z}\theta$, $\theta\in \bold{R}
\setminus\bold{Q}$, oriented by their embedding into $\bold{R}$,
and whose morphisms are strict isomorphisms, that is,
multiplications by a positive number identifying two pseudolattices.
According to (1.2), such isomorphisms 
$L_{\theta'}\to L_{\theta}$ are represented by
matrices 
$$
g=\left(\matrix a&b\\c&d\endmatrix\right)\in GL\,(2,\bold{Z})
$$
such that 
$$
 \theta^{\prime} = g\theta =\frac{a\theta +b}{c\theta +d},\quad
c\theta +d > 0.
\eqno(1.15)
$$ 

\smallskip

Denote by $\Cal{Q}\Cal{T}_{iso}$ the category whose objects
are irrational toric algebras $A_{\theta}$ and whose morphisms
are bimodules inducing Morita equivalences.

\smallskip

Given $\theta',\,\theta ,$ and $g$ satisfying (1.15),
construct an $(A_{\theta'},A_{\theta})$--bimodule 
${}_{\theta'}E_{\theta}(g^{-1})$ (notice the inversion $g^{-1}$)
by the following
prescription. The smooth part of ${}_{\theta'}E_{\theta}(g^{-1})$
consists of functions $f(x,\mu )$ in the Schwartz's space
$\Cal{S} (\bold{R}\times\bold{Z}_c)$. The generators
$U,V$ of $A_{\theta}$ act upon these functions from the right as follows:
$$
(fU)(x,\mu )=f(x-\frac{c\theta +d}{c},\mu -1),
$$
$$
(fV)(x,\mu )=e^{2\pi i (x-\mu d/c)}f(x,\mu ).
$$
The generators $U',V'$ of $A_{\theta'}$ act from the left:
$$
(U'f)(x,\mu )=f(x-\frac{1}{c},\mu -a),
$$
$$
(V'f)(x,\mu )=\roman{exp}\,\left[2\pi i \left(\frac{x}{c\theta +d}
-\frac{\mu}{c}\right)\right]\,f(x,\mu ).
$$
To become a bimodule over the respective $C^*$ algebras, the
Schwartz space must be appropriately completed, cf. Theorem 3.4.1
below. 

\smallskip

\proclaim{\quad 1.7.3. Theorem} The map $E:\, \Cal{P}\Cal{L}_{iso}\to \Cal{Q}\Cal{T}_{iso} $ defined on objects by
$L_{\theta}\mapsto A_{\theta}$  and sending the isomorphism
(1.15) to the Morita isomorphism $[{}_{\theta'}E_{\theta}(g^{-1})]$ 
is a well defined functor. The composition
$K\circ E$ is isomorphic to the identical functor on $\Cal{P}\Cal{L}_{iso}$. 
\endproclaim

\medskip

This theorem rephrases the main result of [DiSch], section 2,
which in our notation establishes an explicit isomorphism
of bimodules
$$
s_{g,h}:\, {}_{\theta'}E_{\theta}(g^{-1})\otimes_{A_{\theta}}
{}_{\theta}E_{\theta''}(h^{-1}) \to
{}_{\theta'}E_{\theta''}((gh)^{-1})
$$
and thus shows that $E$ is multiplicative on isomorphisms 
of pseudolattices. Here $\theta=g^{-1}(\theta' )$ as above,
and $\theta''=h^{-1}(\theta )$ 
so that $\theta''=(gh)^{-1}(\theta').$  

\smallskip

This isomorphism is constructed in [DiSch] in the smooth 
setting. According to [Co2], extension of rings
induces a bijection between the set of isomorphism classes of 
projective modules of finite type over $A_{\theta}$ and 
over its smooth subring respectively. Moreover, the trace (dimension) of
${}_{\theta'}E_{\theta}(g^{-1})$ as a right $A_{\theta}$--module
equals $|c\theta +d|$ ([Co2], Theorem 7): since ${}_{\theta'}E_{\theta}(g^{-1})$
is not given as the image of a projection, Connes develops
differential geometric methods for calculating this trace.
 As an exercise,
the reader can check that the dimension of the tensor product
indeed equals the product of dimensions of factors.
\smallskip

Notice in conclusion that our version of Morita category
using isomorphism classes of bimodules as morphisms
is a truncation of a finer notion which treats
bimodules as functors and leads to the notion
of Morita 2--category. A refinement of the 
Dieng--Schwarz's result in this direction requires
an explicitation of the associativity isomorphism 
connecting $s_{gh,k}\circ (s_{g,h}\otimes \roman{id})$
to $s_{g,hk}\circ (\roman{id}\otimes s_{h,k})$
which replaces the straightforward associativity of the triple
multiplication of morphisms in 1--categories.
This looks like a nice exercise.

\medskip

{\bf 1.8. Quantum tori as ``limits'' of elliptic curves.}
Reading parallelly subsections 0.1 and 1.1,
we see that pseudolattices are in a very precise sense
limits of lattices, at least, if one forgets
orientation; or else one can add orientation 
to the definition of a lattice, as the choice
of a generator of $\wedge^2 (\Lambda ).$

\smallskip

Passing to the isomorphism classes of lattices/pseudolattices
does not seem to change this impression:
compare (0.3) and (1.3).

\smallskip

Comparison of the relevant geometric categories
suggests that two--dimensional quantum tori
can be thus considered as limits of elliptic curves.
More specifically, take a family
of Jacobi parametrized curves $E_{\tau}=\bold{C}/(e^{2\pi i\tau})$
with $\roman{Im}\,\tau > 0$ and $\tau \to \theta \in \bold{R}$.
It is then natural to imagine $T_{\theta}$ as a limit
of $E_{\tau}.$ 

\smallskip

Fixing a Jacobi uniformization of an elliptic curve
(or abelian variety of any dimension) as 
a part of its structure is necessary, for example, 
in problems connected with mirror symmetry.
In such contexts our intuition seemingly provides
a sound picture (cf. a similar discussion in [So],
pp. 100, 113--114).

\smallskip

However, limitations
of this viewpoint become quite apparent if one
has no reason to keep a Jacobi uniformization
as a part of the structure, and is interested only
in the isomorphim classes of elliptic curves, perhaps
somewhat rigidified by a choice of a level structure.

\smallskip

In this case one must contemplate the dynamics of the limiting process
not on the closed upper half--plane but on a relevant modular
curve $X$. Letting $\tau$ tend to $\theta$ along a geodesic,
we get a parametrized real curve on $X$ which, when $\theta$
is irrational, does not tend to any limiting point.
The following lemma shows what can happen.

\medskip

\proclaim{\quad 1.8.1. Lemma} (a) Let $\theta$ be a real quadratic
irrationality, $\theta^{\prime}$ its conjugate.
Consider the oriented geodesic in $H$ joining $\theta^{\prime}$
to $\theta$. The image of this geodesic on any modular curve
$X$ is supported by
a closed loop, which we denote $(\theta^{\prime},\theta )_X$.

\smallskip

(b) Let $\theta$ be as above, and let $\tau$ tend to $\theta$
along an arbitrary geodesic. Then the image of this
geodesic on $X$ has $(\theta^{\prime},\theta )_X$
as a limit cycle (in positive time).

\smallskip

(c) Each closed geodesic on $X$ is the support 
of a closed loop $(\theta^{\prime},\theta )_X$. The union
of them is dense in $X$. It is a strange attractor  
for the geodesic flow in the following sense.
Having chosen a sequence of loops $(\theta^{\prime}_i,\theta_i )_X$,
a sequence of integers $n_i\ge 1$, and a sequence of real
numbers $\epsilon_i>0$, $i=1,2, \dots $, one can find
an oriented geodesic winding $\ge n_i$ times in the 
$\epsilon_i$--neighborhood of $(\theta^{\prime}_i,\theta_i )_X$ for each
$i$, before jumping to the next loop.
\endproclaim

\smallskip

{\bf Proof.} We will
only sketch a couple of arguments. 

\smallskip

For (a), notice that $\theta^{\prime}$ and $\theta$
are respectively the attracting and the repelling
points of a hyperbolic fractional linear transformation
$g\in SL(2,\bold{Z})$. This transformation maps into
itself the whole geodesic joining $\theta^{\prime}$ to $\theta$ 
and acts upon it as a shift by the distance $\roman{log}\,\varepsilon$
where $\varepsilon >1$ is a unit in the quadratic field
generated by $\theta$ (cf. formula (1.16) below). 
If $X=\Gamma \subset H$, where $\Gamma$
is a subgroup of finite index of the modular group,
then $g^n\in \Gamma$ for an appropriate $n\ge 1$. Therefore
the geodesic in question will close to a loop on $X$.

\smallskip

The distance between two geodesics tending to the same
$\theta$ in $H$ tends to zero; this shows (b).

\smallskip

Finally, (c) is based upon an elementary argument
involving continued fractions and diophantine approximations. 
The Lemma is proved.

\medskip

Now let us imagine that we have constructed a certain object
$R(E_\tau )$ depending on the isomorphism class of $E_\tau$
(perhaps, with rigidity). This object can be a number,
a function of the lattice, a linear space, a category ... Suppose also
that we have constructed a similar object $\Cal{R} (T_{\theta})$
depending on the isomorphism class of $T_{\theta}$,
and that we want to make sense of the intuitive notion
that $\Cal{R} (T_{\theta})$ is ``a limit of $R(E_\tau )$.''
Since in the most interesting for us case (a) of the Lemma 1.8.1
$E_{\tau}$ keeps rotating around the same loop,
there are  two natural possibilities:

\medskip

(i) {\it The object $R(E_\tau )$ actually ``does not depend on
$\tau$'', and $\Cal{R} (T_{\theta})$ is its constant
value. Here independence generally means a canonical identification
of different $R(E_\tau )$, e.g. via a version of  flat
connection defined along the loop.} 

\smallskip

(ii)  {\it The object $R(E_\tau )$ does depend on
$\tau$, and $\Cal{R} (T_{\theta})$ is obtained by a 
kind of integrating or averaging various
$R(E_\tau )$ along the loop.}

\medskip

The second case looks more interesting, however, it
is not immediately
obvious that such objects occur in nature.
Remarkably, they do, and precisely in the context
of real multiplication and Stark's conjecture.
In fact, this is how we will interpret the
beautiful old calculational tricks  due to Hecke:
see [He1], [He2], [Her], [Z]. See also [Dar]
for a similar observation related to what Darmon
calls Stark--Heegner points of elliptic curves.

\smallskip

In this section we will
only explain the geometric meaning of Hecke's substitution,
whereas the (slightly generalized) calculation itself will be treated
in the next section.

\medskip

{\bf 1.8.2. Hecke's lift of closed
geodesics to the space of lattices.} Let $K\subset \bold{R}$
be a real quadratic subfield of $\bold{R}$ and
$L\subset K$ an RM pseudolattice. From now on, we denote
by $l\mapsto l^{\prime}$ the nontrivial element
of the Galois group of $K/\bold{Q}$. 

\smallskip

For any real $t$, consider the following subset of $\bold{C}$:
$$
\Lambda_t = \Lambda_t (L) :=
\{\lambda_t=\lambda_t(l):=le^{t/2}+il^{\prime}e^{-t/2}\,|\, l\in L\}
\eqno(1.16)
$$

\smallskip

\proclaim{\quad Lemma 1.8.3} (a) $\Lambda_t(L)$ is a lattice.

\smallskip

(b) Any isomorphism $a:\,L_1\to L$ in the narrow sense
induces  isomorphisms $\Lambda_t(L_1)\to \Lambda_{t+c}(L)$
where $c$ is a constant depending only on $a$
and $t$ is arbitrary.
\smallskip

(c) The image of the curve $\{\Lambda_t\,|\,t\in\bold{R}\}$
on the modular curve (0.3) (or any modular curve) is a closed
geodesic. The affine coordinate $t$ along this curve
is the geodesic length.
\endproclaim

\smallskip

{\bf Proof.} (a) is evident; moreover, if $l_1,l_2$
form a basis of $L$, then $\lambda_t(l_1),\lambda_t(l_2)$
form a basis of $\Lambda_t$.

\smallskip

For (b), consider an isomorphism $L\mapsto L_1:\,
l\mapsto al,\,a\in K,\, aa^{\prime}>0.$ It induces a map
$\Lambda_t(L)\to \Lambda_t(L_1)$:
$$
\lambda_t(l)\mapsto al\,e^{t/2}+ia^{\prime}l^{\prime}e^{-t/2}=
$$
$$
\sqrt{aa^{\prime}}\,\left(\sqrt{\frac{a}{a^{\prime}}}\,le^{t/2}
+\sqrt{\frac{a^{\prime}}{a}}\,l^{\prime}e^{-t/2}\right)=
\sqrt{aa^{\prime}}\,\lambda_{t+\roman{log}\,\frac{a}{a^{\prime}}}(l).
\eqno(1.17)
$$
This produces an isomorphism of $\Lambda_t(L_1)$ with
$\Lambda_{t+\roman{log}\,\frac{a}{a^{\prime}}}(L).$

\smallskip

For (c), it suffices to consider pseudolattices $L$
generated by $1$ and $\theta\in K$ with $\theta^{\prime}>\theta .$
Then $\Lambda_t(L)$ is generated by
$e^{t/2}+ie^{-t/2}$ and $\theta\,e^{t/2}+i\theta^{\prime}\,e^{-t/2}$,
and hence isomorphic to the lattice generated
by $1$ and
$$
\tau_t:=\frac{\theta\,e^{t/2}+i\theta^{\prime}\,e^{-t/2}}{e^{t/2}+ie^{-t/2}}=
\frac{\theta\,e^{t}+\theta^{\prime}\,e^{-t}}{e^{t}+e^{-t}}
+i\,\frac{\theta^{\prime}-\theta}{e^{t}+e^{-t}}.
\eqno(1.18)
$$
A straightforward computation shows that
$$
\left|\,\tau_t-\frac{\theta +\theta^{\prime}}{2}\,\right|^2=
\left( \frac{\theta^{\prime} -\theta}{2}\right)^2.
$$
Hence $\tau_t$ runs over a semicircle in the upper half
plane connecting $\theta^{\prime}$ to $\theta$.
A further calculation shows that the geodesic length element
$\dfrac{|\,d\tau|}{\roman{Im}\,\tau}$ restricted to
this semicircle coincides with $dt$. The normalization of $t$ has a 
simple geometric meaning: $t=0$ is the upper point
of the geodesic semicircle.

\bigskip


\centerline{\bf \S 2. Stark's numbers and theta functions}

\smallskip

\centerline{\bf for real quadratic fields}

\medskip

{\bf 2.1. Stark's numbers at $s=0$.} In this section we fix
a real quadratic subfield $K\subset \bold{R}$. Denote
by $l\mapsto l^{\prime}$  the action of the nontrivial
element of the Galois group of $K$, and by $O_K$ the ring
of integers of $K$, and put $N(l)=ll^{\prime}.$

\smallskip

Let $L$ be an arbitrary integral ideal of $K$ which,
together with its embedding in $\bold{R}$ and the induced
ordering, will be considered as a pseudolattice.

\smallskip

Choose also an $l_0\in O_K$ so that the pair $(L,l_0)$
satisfies the following restrictions:

\medskip

(i) {\it The ideals $\frak{b}:=(L,l_0)$ and $\frak{a}_0:=(l_0)\frak{b}^{-1}$
are coprime with $\frak{f}:=L\frak{b}^{-1}$.}

\smallskip

(ii) {\it Let $\varepsilon$ be a unit of $K$ such that
$\varepsilon\equiv 1\,\roman{mod}\,\frak{f}.$ Then $\varepsilon^{\prime}>0.$}

\medskip

Put now
$$
\zeta (L,l_0,s):= \roman{sgn}\,l_0^{\prime}\,N(\frak{b})^s\sum^{(u)}_{l\in L}
\frac{\roman{sgn}\,(l_0+l)^{\prime}}{|N(l_0+l)|^s}
\eqno(2.1)
$$
where $(u)$ at the summation sign means that one should
take one representative from each coset $(l_0+l)\varepsilon$
where $\varepsilon$ runs over all units $\equiv\,1\,\roman{mod}\,\frak{f}$.
Notice that $(l_0+L)\varepsilon =l_0+L$ precisely for such units.

\smallskip

With this conventions, our $\zeta (L,l_0,s)$ is exactly
Stark's function denoted $\zeta(s,\frak{c})$
on the page 65 of [St1]: our $\frak{a}_0,\frak{b},\frak{f}$
have the same meaning in [St1], and our $l_0$ is Stark's $\gamma$.
The meaning of Stark's $\frak{c}$ is explained below.

\smallskip

The Stark number of $(L,l_0)$ is defined as
$$
S_0(L,l_0):=e^{\zeta^{\prime}(L,l_0,0)}
\eqno(2.2)
$$
(cf. the general discussion in 0.6).

\smallskip

The simplest examples correspond to the cases when 
$(L,l_0)=(1)$, $\frak{f}=L$, in particular, $l_0=1$.

\smallskip

Notice that pseudolattices which are integral ideals
have conductor $f=1$ in the sense of Lemma 1.1.1.

\medskip

{\bf 2.2. Stark's conjecture for real quadratic fields.} 
In [St1], Stark conjectures
that $S_0(L,l_0)$ are algebraic units generating
abelian extensions of $K$. To be more precise,
let us first describe an abelian extension $M/K$
associated with $(L,l_0)$ using the classical
language of class field theory. (Our $M$ is Stark's $K$,
whereas our $K$ corresponds to Stark's $k$.)

\smallskip

In 2.1 above we constructed, starting with $(L,l_0)$,
the ideals $\frak{f}$ and $\frak{b}$ in $O_K$.
Let $I(\frak{f})$ be the group of fractional ideals of $K$
generated by the prime ideals of $K$ not dividing $\frak{f}$,
and $S(\frak{f})$ be its subgroup called the
principal ray class modulo $\frak{f}$. Then Artin's 
reciprocity map identifies 
$G(\frak{f}):=I(\frak{f})/S(\frak{f})$ with the Galois group of $M/K$.

\smallskip

Consider all pairs $(L,l_0)$ as above with fixed $\frak{f}.$
It is not difficult to establish that on this set,
$S_0(L,l_0)$ in fact depends only on the class 
$\frak{c}$ of
$\frak{a}_0=(l_0)\frak{b}^{-1}$ in $G(\frak{f})$. Denote the respective number $E(\frak{c}).$

\medskip

{\bf 2.2.1. Conjecture.} {\it The numbers $E(\frak{c})$ are units
belonging to $M$ and generating $M$ over $K$.
If the Artin isomorphism associates with $\frak{c}$ an
automorphism $\sigma$, we have $E(1)^{\sigma}=E(\frak{c}).$}

\smallskip

(We reproduced here the most optimistic form of the Conjecture 1
on page 65 of [St1] involving $m=1$ and Artin's  reciprocity map).

\medskip

{\bf 2.3. Hecke's formulas.}  In this subsection we will
work out Hecke's approach to the computation of sums
of the type (2.1), cf. [He2]. It starts with a
Mellin transform so that instead of Dirichlet series (2.1)
we will be dealing with a version of theta--functions
for real quadratic fields. We start with introducing
a class of such theta functions more general than
strictly needed for dealing with  (2.1) (and more general
than Hecke's one). 

\medskip

{\bf 2.4. Theta functions of pseudolattices.}
Let $K\subset \bold{R}$ be as in 2.1. We choose and fix the
following data: a pseudolattice $L\subset K$,
two numbers $l_0,m_0\in K$ and a number
$\eta =\eta_0+i\eta_1\in\bold{C}.$ A complex variable $v$ will take values in
the upper half plane; $\sqrt{-iv}$ is the branch which is positive
on the upper part of the imaginary axis.

\smallskip

Finally, choose an infinite cyclic group $U$ of totally
positive units
in $K$ such that the following conditions hold:

\smallskip

(a) $u(l_0+L)= l_0+L$ for all
$u\in U.$

\smallskip

(b) $\roman{tr}\,ulm_0\equiv\roman{tr}\,lm_0\,\roman{mod}\,\bold{Z}$,
$\roman{tr}\,ul_0m_0\equiv\roman{tr}\,l_0m_0\,\roman{mod}\,2\bold{Z}$
for all $l\in L$, $u\in U,$ where $\roman{tr}:=\roman{tr}_{K/\bold{Q}}.$

\smallskip

Put now
$$
\Theta_{L,\eta}^U
\thickfracwithdelims[]\thickness0{l_0}{m_0}(v):=
$$
$$
\sum_{l_0+l\,\roman{mod}\,U}
(\eta_0\,\roman{sgn}\,(l_0^{\prime}+l^{\prime}) +
\eta_1\,\roman{sgn}\,(l_0+l))\,
e^{2\pi i\,v|(l_0+l)(l_0^{\prime}+l^{\prime})|}
e^{-2\pi i\,\roman{tr}\,lm_0}
e^{-\pi i\,\roman{tr}\,l_0m_0} . 
\eqno(2.3)
$$
Notation $l_0+l\,\roman{mod}\,U$ means that we sum over
a system of representatives of orbits of $U$
acting upon $l_0+L$.

\smallskip

Notice that such $U$ always exists, and that if we choose
a smaller subgroup $V\subset U$, then
$$
\Theta_{L,\eta}^V
\thickfracwithdelims[]\thickness0{l_0}{m_0}(v)=
[U:V]\,\Theta_{L,\eta}^U
\thickfracwithdelims[]\thickness0{l_0}{m_0}(v).
$$

\smallskip

In order to relate these thetas to Stark's numbers,
consider the function
$$
\Theta_{L,1}^U
\thickfracwithdelims[]\thickness0{l_0}{0}(v)=
\sum_{l_0+l\,\roman{mod}\,U}
\roman{sgn}\,(l_0^{\prime}+l^{\prime})\,
e^{2\pi i\,v|(l_0+l)(l_0^{\prime}+l^{\prime})|} .
\eqno(2.4)
$$
Then we have
$$
\sum_{l_0+l\roman{mod}\,U}
\frac{\roman{sgn}\,(l_0^{\prime}+l^{\prime})}{|N(l_0+l)|^s}=
\frac{(2\pi)^s}{\Gamma (s)}\,
\int_0^{i\infty}(-iv)^{s}
\Theta_{L,1}^U
\thickfracwithdelims[]\thickness0{l_0}{0}(v)\,\frac{dv}{v} .
\eqno(2.5)
$$
We will now show that these RM thetas can be obtained
by averaging some theta constants (related to the complex lattices)
along the closed geodesics described in 1.8 above.

\medskip

{\bf 2.5. Theta constants along geodesics.} Starting with the same data
as in 2.4, we introduce first of all a family
of lattices $\Lambda_t=\Lambda_t(L)$ defined by (1.16).
From $l_0$ which was used  to shift $L$, we will produce
a shift of $\Lambda_t$:
$$
\lambda_{0,t}:=l_0\,e^{t/2}+il_0^{\prime}\,e^{-t/2}.
$$
The number $m_0$  determines a character of $L$
appearing in (2.3): $l\mapsto e^{-2\pi i\,\roman{tr}\,lm_0}$. 
Similarly, we will produce a character of $\Lambda_t$
from
$$
\mu_{0,t}:=m_0\,e^{t/2}+im_0^{\prime}\,e^{-t/2}
$$
by using the scalar product on $\bold{C}$
$$
(x\cdot y)=\roman{Im}\, xy= x_0y_1+x_1y_0 
\eqno(2.6)
$$
where $x=x_0+ix_1,\,y=y_0+iy_1.$ Since $l_0,m_0\in L\otimes\bold{Q}$,
we have similarly $\lambda_{0,t}, \mu_{0,t}\in\Lambda_t\otimes\bold{Q}$.
Omitting $t$ for brevity, we put:
$$
\theta_{\Lambda,\eta}
\thickfracwithdelims[]\thickness0{\lambda_0}{\mu_0}\,(v):=
\sum_{\lambda\in\Lambda}
((\lambda_{0}+\lambda)\cdot\eta )\,e^{\pi i v|\lambda_{0}+\lambda|^2}
e^{-2\pi i (\lambda\cdot \mu_{0})-\pi i (\lambda_{0}\cdot
\mu_{0})} .
\eqno(2.7)
$$
The two types of thetas are related by Hecke's averaging formula:

\medskip

\proclaim{\quad 2.6. Proposition} We have
$$
\Theta_{L,\eta}^U
\thickfracwithdelims[]\thickness0{l_0}{m_0}(v)=
\sqrt{-iv}\,
\int_{-\roman{log}\,\varepsilon}^{\roman{log}\,\varepsilon} 
\theta_{\Lambda_t,\eta}
\thickfracwithdelims[]\thickness0{\lambda_{0,t}}{\mu_{0,t}}\,(v)\, dt 
\eqno(2.8)
$$
where $\varepsilon >1$
is a generator of $U$. 
\endproclaim

{\bf Proof.} The following formulas are valid
for $\roman{Im}\,v>0$:
$$
e^{2\pi i\,v|mm^{\prime}|} =
\sqrt{-iv}\,|m^{\prime}|\int_{-\infty}^{\infty}
e^{-t/2}e^{\pi iv(m^2e^{t}+ m^{\prime 2}e^{-t})}\,dt =
$$
$$
\sqrt{-iv}\,|m|\int_{-\infty}^{\infty}
e^{t/2}e^{\pi iv(m^2e^{t}+ m^{\prime 2}e^{-t})}\,dt
\eqno(2.9)
$$
(see e.~g. [La], pp. 270--271).  In the rhs of (2.3), 
replace the first exponent
by its integral versions (2.9), using the first version
at $\eta_0$ and the second at $\eta_1$. We get:
$$
\Theta_{L,\eta}^U
\thickfracwithdelims[]\thickness0{l_0}{m_0}(v)=
$$
$$
\sqrt{-iv}\int_{-\infty}^{\infty} 
\sum_{l_0+l\,\roman{mod}\,U}
(\eta_0\,(l_0^{\prime}+l^{\prime})\,e^{-t/2} +
\eta_1\,(l_0+l)\,e^{t/2})\,\times
$$
$$
e^{\pi iv((l_0+l)^2e^{t}+ (l^{\prime}_0+l^{\prime})^2e^{-t})}
e^{-2\pi i\roman{tr}\,lm_0}
e^{-\pi i\roman{tr}\,l_0m_0}\,dt\,.
\eqno(2.10)  
$$
In view of (1.16) and (2.6) we have
$$
\eta_0\,(l_0^{\prime}+l^{\prime})\,e^{-t/2} +
\eta_1\,(l_0+l)\,e^{t/2} =
((\lambda_{0,t}+\lambda_{t})\cdot\eta ),
$$
$$
(l_0+l)^2e^{t}+ (l^{\prime}_0+l^{\prime})^2e^{-t}=
|\lambda_{0,t}+\lambda_t|^2,
$$
and similarly
$$
\roman{tr}\,lm_0=(\lambda_t\cdot\mu_{0,t}),\quad
\roman{tr}\,l_{0}m_{0}=(\lambda_{0,t}\cdot\mu_{0,t}) .
$$
Inserting this into (2.10), we obtain
$$
\sqrt{-iv}\,\int_{-\infty}^{\infty} \,dt
\sum_{l_0+l\,\roman{mod}\,U}
((\lambda_{0,t}+\lambda_t)\cdot\eta)
e^{\pi i\, v|\lambda_{0,t}+\lambda_t|^2}\,
e^{-2\pi i\, (\lambda_t\cdot\mu_{0,t})}
e^{-\pi i\, (\lambda_{0,t}\cdot\mu_{0,t})} .
\eqno(2.11)  
$$
Replacing $l_0+l$ by $\varepsilon (l_0+l)$ is equivalent
to replacing $t$ by $t+2\,\roman{log}\,\varepsilon$.
Hence finally the right hand side of (2.11) can be rewritten
as
$$
\sqrt{-iv}\,
\int_{-\roman{log}\,\varepsilon}^{\roman{log}\,\varepsilon} dt 
\sum_{\lambda_t\in\Lambda_t}
((\lambda_{0,t}+\lambda_t)\cdot\eta )\,e^{\pi i\, v|\lambda_{0,t}+\lambda_t|^2}
e^{-2\pi i\, (\lambda_t\cdot \mu_{0,t})-\pi i\, (\lambda_{0,t}\cdot
\mu_{0,t})}
\eqno(2.12)
$$
which is the same as (2.8).

\smallskip

We will now apply Poisson formula in order to derive
functional equations for Hecke's thetas.

\medskip

{\bf 2.7. Poisson formula.} Let $V$ be a real vector space,
$\widehat{V}$ its dual. We will denote by $(x\cdot y)\in\bold{R}$ the scalar
product of $x\in V$ and $y\in \widehat{V}.$
Choose a lattice (discrete subgroup of finite covolume) $\Lambda\subset
V$ and put
$$
\Lambda^!:=\{\mu \in \widehat{V}\,|\,\forall \lambda \in \Lambda,\,
(\lambda\cdot\mu )\in\bold{Z}\} .
\eqno(2.13)
$$
Choose also a Haar measure $dx$ on $V$ and define the Fourier transform
of a Schwarz function $f$ on $V$ by
$$
\widehat{f}(y):=\int_V f(x)\,e^{-2\pi i (x\cdot y)} dx . 
\eqno(2.14)
$$
If $f(x)$ in this formula is replaced by
$f(x+x_0)\,e^{-2\pi i (x\cdot y_0)-\pi i (x_0\cdot y_0)}$
for some $x_0\in V,\,y_0\in \widehat{V}$,
its Fourier transform $\widehat{f}(y)$ gets replaced
by $\widehat{f}(y+y_0)\,e^{2\pi i (x_0\cdot y)+\pi i (x_0\cdot y_0)}.$

\smallskip

The Poisson formula reads
$$
\sum_{\lambda\in\Lambda} f(\lambda )=\frac{1}{\int_{V/\Lambda} dx}\,
\sum_{\mu\in\Lambda^!} \widehat{f}(\mu ) ,
\eqno(2.15)
$$
and for shifted functions as above
$$
\sum_{\lambda\in\Lambda} f(\lambda_0+\lambda)\,e^{-2\pi i (\lambda\cdot \mu_0)-\pi i (\lambda_0\cdot \mu_0)} =\frac{1}{\int_{V/\Lambda} dx}\,
\sum_{\mu\in\Lambda^!} \widehat{f}(\mu_0+\mu )\,e^{2\pi i (\lambda_0\cdot \mu )+\pi i (\lambda_0\cdot \mu_0) } .
\eqno(2.16)
$$
\medskip
{\bf 2.8. Functional equations for $\theta$ and $\Theta$.}
In order to transform (2.12) using the Poisson formula,
we put
$$
V=\bold{C}=\{x_0+ix_1\},\ \widehat{V}= \bold{C}=\{y_0+iy_1\},\
\eqno(2.17)
$$
and take (2.6) for the scalar product.

\smallskip

\proclaim{\quad 2.8.1. Lemma} Let the lattice $\Lambda_t\subset \bold{C}$
be given by (1.15). Then the dual lattice $\Lambda_t^!$ with respect
to the pairing (2.6) has the similar structure
$$
\Lambda_t^! = \Lambda_t (M) :=
\{me^{t/2}+im^{\prime}e^{-t/2}\,|\, m\in M\} 
\eqno(2.18)
$$
where we denoted by $M=L^?$ the pseudolattice
$$
M:=\{ m\in K\,|\,\forall l\in L,\,\roman{tr}_{K/\bold{Q}} (l^{\prime}m)
\in\bold{Z}.\}.
$$
\endproclaim
\smallskip

{\bf Proof.} Denote by $\Gamma$ the lattice (2.18).
For any $\lambda =  le^{t/2}+il^{\prime}e^{-t/2}\in 
\Lambda_t$ and $\mu = me^{t/2}+im^{\prime}e^{-t/2}\in \Gamma$
we have
$$
(\lambda\cdot\mu )= \roman{Im}\,\lambda\mu=lm^{\prime}+l^{\prime}m=
\roman{tr}_{K/\bold{Q}}(lm^{\prime}).
\eqno(2.19)
$$
Therefore this scalar product lies in $\bold{Z}$ if
$m\in M$ so that $\Gamma\subset \Lambda_t^!$.
Clearly, then, $\Gamma$ must be commensurable with $\Lambda_t^!$,
so that the right hand side of (2.19) can be used
for computing  $(\lambda\cdot\mu )$ on the whole $\Lambda_t^!$.
This finishes the proof.

\smallskip
For example,  $O_K^? =\frak{d}^{-1}$ where
$\frak{d}$ is the different. In fact,
this is the standard definition of the different. 

\smallskip

Now let $l_1,l_2$ be two generators of the pseudolattice $L$.
Put
$$
\Delta (L):= |l_1l_2^{\prime}-l_1^{\prime}l_2| .
\eqno(2.20)
$$
Clearly, this number does not depend on the choice of
generators.

\medskip

\medskip

\proclaim{\quad 2.8.2. Lemma} Let the Haar measure on $V$ be
$dx=dx_0\,dx_1.$ Choose generators $l_1,l_2$ of $L$. Then
$$
\int_{V/\Lambda_t}dx = \Delta (L) .
\eqno(2.21)
$$
\endproclaim
\smallskip

{\bf Proof.} If $\Lambda_t$ is generated by $\omega_1,\omega_2$,
then the volume (2.21) equals 
$$
|\roman{Re}\,\omega_1 \,\roman{Im}\,\omega_2 -
\roman{Re}\,\omega_2 \,\roman{Im}\,\omega_1 |.
$$
Taking
$$
\omega_1=l_1e^{t/2}+il_1^{\prime}e^{-t/2},\quad
\omega_2=l_2e^{t/2}+il_2^{\prime}e^{-t/2},
$$
we get (2.20).

\medskip

\proclaim{\quad 2.8.3. Lemma} The Fourier transform of
$$
f_{v,\eta}(x):= (x\cdot\eta )\,e^{\pi i v|x|^2},\ \eta =\eta_0+i\eta_1
\eqno(2.22)
$$
equals
$$
{g}_{v,\eta}(y):=\frac{i}{v^2}\,(y\cdot i{\bar\eta})\,e^{-\frac{\pi i}{v}|y|^2}
\eqno(2.23)
$$
\endproclaim
\smallskip

{\bf Proof.} Putting $w=-iv$ we have
$$
f_{v,\eta}(x)=(x_0\eta_1+x_1\eta_0)\,e^{-\pi w(x_0^2+x_1^2)},
$$
so that its Fourier transform by (2.13) and (2.14) is
$$
\eta_1 \int_{-\infty}^{\infty} e^{-\pi wx_0^2}\,e^{-2\pi i x_0y_1}\, x_0\,dx_0\,\cdot\,
\int_{-\infty}^{\infty} e^{-\pi wx_1^2}\,e^{-2\pi i x_1y_0}\,dx_1 +
$$
$$
\eta_0 \int_{-\infty}^{\infty} e^{-\pi wx_0^2}\,e^{-2\pi i x_0y_1}\, dx_0\,\cdot\,
\int_{-\infty}^{\infty} e^{-\pi wx_1^2}\,e^{-2\pi i x_1y_0}\,x_1\,dx_1 =
$$
$$
(\eta_0y_0+\eta_1y_1)\,\frac{1}{iw^2}\,e^{-\pi\frac{y_0^2+y_1^2}{w}} .
$$
This is (2.23).

\medskip

{\bf 2.8.4. A functional equation for $\theta$.} Let us now 
write (2.16) for $f=f_{v,\eta}$ and $\Lambda_t$:
$$
\sum_{\lambda\in\Lambda_t}
((\lambda_{0,t}+\lambda)\cdot\eta )\,e^{\pi i v|\lambda_{0,t}+\lambda|^2}
e^{-2\pi i (\lambda\cdot \mu_{0,t})-\pi i (\lambda_{0,t}\cdot
\mu_{0,t})}=
$$
$$
\frac{i}{\Delta (L)\,v^2}\,\sum_{\mu\in\Lambda_t^!}
((\mu_{0,t}+\mu)\cdot i\bar{\eta} )\,e^{-\frac{\pi i} {v} |\mu_0+\mu |^2}
e^{2\pi i (\lambda_{0,t}\cdot \mu )+\pi i (\lambda_{0,t}\cdot
\mu_{0,t})} .
$$ 
In the notation (2.7) this means:
$$
\theta_{\Lambda_t,\eta}
\thickfracwithdelims[]\thickness0{\lambda_{0,t}}{\mu_{0,t}}\,(v)=
\frac{i}{\Delta (L)\,v^2}\,
\theta_{\Lambda_t^!,i\bar{\eta}}
\thickfracwithdelims[]\thickness0{\mu_{0,t}}{-\lambda_{0,t}}\,
\left(-\frac{1}{v}\right) .
\eqno(2.24)
$$
\smallskip

We now can establish a functional equation
for $\Theta^U$ as well:

\medskip

\proclaim{\quad 2.9. Proposition} We have
$$
\Theta_{L,\eta}^U
\thickfracwithdelims[]\thickness0{l_0}{m_0}(v)=
\frac{1}{\Delta (L)\,v}\,
\Theta_{L^?,i\bar{\eta}}^U
\thickfracwithdelims[]\thickness0{m_0}{-l_0}\left(-\frac{1}{v}\right) .
\eqno(2.25)
$$
\endproclaim
\smallskip

{\bf Proof.} This is a straightforward consequence of (2.8)
and (2.24).

\bigskip


\centerline{\bf \S 3. Heisenberg groups, modules over quantum tori,}
\smallskip

\centerline{\bf and theta functions}

\medskip

{\bf 3.0. Introduction.} Most of the constructions of this
section are explained for the case of tori
of arbitrary dimension. In 3.1--3.5 we remind
to the reader the approach to the classical theta functions
based upon the theory of Heisenberg groups. 
We closely
follow Mumford's presentation in [Mu3], \S 1 and \S 2,
which ideally suits our goals. The reader can
find missing proofs there.

\smallskip

Quantum tori and their representations appear very 
naturally, when one restricts the basic Heisenberg
representation to a lattice. This leads naturally to
the emergence of Rieffel's setup as in Lemma 1.5.1, (1.8)--(1.11),
although no explicit projections form a part of the picture.
A way to remedy this and to construct certain projections
starting with theta functions was proposed by F.~Boca in [Bo2].
Generalizing his calculation, we prove the Theorem 3.7,
which introduces in the context of representation theory
of toric algebras {\it quantum thetas} in the sense of [Ma3].
This is the third type of thetas we meet in this paper
(counting $\Theta^U$ and $\theta$ of Section 2 for the first two),
and thanks to Boca's theorem, they can be used
to construct morphisms of quantum tori.

\smallskip

The initial motivation of [Ma3] was to produce
quantized versions of coordinate rings of abelian
varieties, generated by the classical
theta constants, i.~e. the values of theta fuctions
at the toric points of finite order. The way
they appear here gives a partial answer to the question raised
by A.~S.~Schwarz in [Sch2].

\medskip

{\bf 3.1. Heisenberg groups.} 
We start with a locally compact abelian topological group
$\Cal{K}$ and denote its character group $\widehat{\Cal{K}}=
\roman{Hom}\,(\Cal{K},\bold{C}^*_1),$ $\bold{C}^*_1=
\{ z\in \bold{C}\,|\,|z|=1\}.$ We also choose
a skew--symmetric pairing $\epsilon :\, \Cal{K}\times \Cal{K}\to 
\bold{C}^*_1$ which is {\it non--degenerate} in the following
sense: it induces an isomorphism $\Cal{K}\to \widehat{\Cal{K}}$,
and $\epsilon (x,x)\equiv 1$. 

\smallskip

Moreover, choose a compatible with $\epsilon$ cocycle $\psi :\,K\times K\to \bold{C}^*_1$:
$$
\psi (x,y) \psi (x+y,z) = \psi (x,y+z) \psi (y,z),
\eqno(3.1)
$$
$$
\epsilon (x,y) = \frac{\psi (x,y)}{\psi (y,x)} .
\eqno(3.2)
$$
The condition (3.1) holds automatically if $\psi$
is a bicharacter. Hence if one can find a skewsymmetric bicharacter
$\epsilon^{1/2}$ which is a square root of $\epsilon$, it can be taken for $\psi$.

\smallskip

Another useful construction starts with $\Cal{K}$ which is
already represented as $K_0\times \widehat{K}_0$ for a
topological group $K_0$. Denoting by
$\langle *,*\rangle :\, K_0\times \widehat{K}_0 \to \bold{C}^*_1$
the canonical pairing, we can simultaneously put
$$
\psi ((x,\widehat{x}), (y,\widehat{y})):= \langle x,\widehat{y}\rangle ,\quad
\epsilon ((x,\widehat{x}), (y,\widehat{y})) =
\frac{\langle x,\widehat{y}\rangle}{\langle y,\widehat{x}\rangle}.
\eqno(3.3)
$$

\smallskip

Having chosen $\Cal{K}$ and $\psi$, we can construct the following objects:

\medskip

(i) {\it The Heisenberg group $\Cal{G}=\Cal{G} (\Cal{K},\psi )$.}

\smallskip

As a set, $\Cal{G}$ is $\bold{C}^*_1\times \Cal{K}$, and the composition
law is given by
$$
(\lambda ,y) (\mu ,z) =(\lambda\mu \psi (y,z) ,y+z).
\eqno(3.4)
$$
The associativity is assured by (3.1). The group comes 
as a central extension
$$
1\to \bold{C}_1^*\to \Cal{G}\to K \to 1.
\eqno(3.5)
$$
If $\Cal{K}$ and $\psi$ split as in  (3.3),
both subgroups $K_0$ and $\widehat{K}_0$ of $\Cal{K}$ come together 
with their lifts to $\Cal{G}$: $x\mapsto (1,x).$

\medskip

(ii) {\it Representations of $\Cal{G}$ on functions on $\Cal{K}$.}

\smallskip
Consider a linear space of complex ``functions'' on $\Cal{K}$
which is stable with respect to all shifts $s_x$, $(s_xf)(y)=f(x+y),$
$x,y\in\Cal{K}$. Here the  word ``functions''
should be understood liberally: completions
of spaces of usual functions and  distributions will do as long as shifts
can be extended in such a way that $s_xs_y=s_{x+y}$, and notation $f(x+y)$ does not imply
that we want literally take values at points.

\smallskip

In this case the formula
$$
(U_{(\lambda ,y)}f)(x) :=\lambda \psi (x,y) f(x+y)
\eqno(3.6)
$$
determines a linear representation of $\Cal{G}$ on this space.

\medskip

{\bf 3.2. Basic unitary representation  and its various models.}
In the notation above, a closed subgroup $K_0\subset \Cal{K}$ is called
{\it isotropic}, if $\epsilon (x,y)=1$ for all $x,y\in K_0$, and
{\it maximal isotropic}, if
$K_0$ is maximal with this property. One can then lift
$K_0$ to $\Cal{G}$, i.e. to find a homomorphism
$K_0\to\Cal{G}:\,x\mapsto (\gamma (x), x).$

\smallskip

Assume that such $K_0$ and $\gamma$ are fixed. Consider
the subspace $\Cal{H}(K_0,\gamma )\subset L_2(\Cal{K})$ consisting of
all functions satisfying the condition
$$
\forall\, y\in K_0, \quad f(x+y)=\gamma (y)^{-1}\psi (y, x)^{-1} f(x).
\eqno(3.7)
$$
Using (3.6), this can be equivalently written as
$$
\forall\, y\in K_0, \quad (U_{(\gamma (y),y)}f)(x)=\epsilon (x,y) f(x).
\eqno(3.8)
$$
A straightforward calculation shows that this space is
invariant with respect to the operators (3.6)
and therefore determines a unitary representation of $\Cal{G}$. 

\smallskip

In the particular case when $\Cal{K}$ is $K_0\times \widehat{K}_0$
and the cocycle is as in (3.3), we can identify 
$\Cal{H}(K_0,\gamma )$ with $L_2(\widehat{K}_0)$ because
(3.8) allows us to reconstruct any function from
its restriction to $\widehat{K}_0$.

\smallskip

This construction plays the central role in the theory
of Heisenberg groups because of the following two
facts:

\medskip

\proclaim{\quad 3.2.1. Theorem} (a) $\Cal{H}(K_0,\gamma )$
is irreducible. 

\smallskip

(b)  Any other unitary
representation of $\Cal{G}$ whose restriction on $\bold{C}_1^*$
is $U_{(\lambda ,0)}=\lambda\,\roman{id}$
is isomorphic to the completed
tensor product of $\Cal{H}(K_0,\gamma )$ and a trivial representation.
In particular, representations $\Cal{H}(K_0,\gamma )$
for different choices of $(K_0,\gamma )$ are all isomorphic.
\endproclaim

\smallskip

The non--degeneracy of $e$ is essentially used in the
proof of this unicity statement. Everything said in 3.1
holds without any non--degeneracy assumption.

\medskip

{\bf 3.3. Heisenberg groups and modules over quantum tori.} Since in this
section we will be dealing with quantum tori of arbitrary
dimension, it is convenient to introduce some
invariant notation. Let $D$ be a free abelian
group of finite rank and $\alpha :\,D\times D\to 
\bold{C}_1^*$ a skewsymmetric pairing.
The $C^*$ algebra $C(D,\alpha )$ of the quantum torus $T(D,\alpha )$
with the character group $D$ and quantization
parameter $\alpha$ is the universal algebra
generated by the family of unitaries $e(h)=e_{D,\alpha}(h),\,h\in D$,
satisfying the relations
$$
e(g) e(h)= \alpha (g,h) e(g+h).
\eqno(3.9)
$$
(Left) modules over such tori can be obtained by the following
construction: choose a Heisenberg group $\Cal{G}(\Cal{K},\psi )$
with a bicharacter cocycle  $\psi$ and compatible $\epsilon$.
Consider  a lattice embedding
$l:\,D\hookrightarrow \Cal{K}$. Denote by $\alpha_D$
the bicharacter on $D$ induced by $\epsilon$. Choose a basic
representation $U$ of $\Cal{G}(\Cal{K},\psi )$ in the space $\Cal{H}$
and define the action of $C(D,\alpha )$ on $\Cal{H}$
by
$$
e_{D,\alpha}(h) f:= U_{(1,l(h))}f
\eqno(3.10)
$$
It turns
out that an appropriate completion
of the subspace of smooth functions is a projective
module (see Rieffel's Theorem 3.4.1 below).

\medskip

{\bf 3.4. Basic representations as toric bimodules.} In the setup
of the last paragraph, assume to shorten notation that
$D$ is a lattice (discrete subgroup with compact quotient)  in $\Cal{K}$ and denote by $D^!$ the 
the maximal orthogonal subgroup:
$$
D^!:= \{x\in \Cal{K}\,|\,\forall\,h\in D, \epsilon (h,x)=1\}.
\eqno(3.11)
$$
Let $\alpha^!$ be the pairing induced by $\epsilon$ on $D^!$.
If $D^!$ is free of finite rank, we get similarly
the representation of $C(D^!,\alpha^!)$ on $\Cal{H}$.
Moreover, operators from $C(D,\alpha)$  and $C(D^!,\alpha^!)$
pairwise commute. Identifying $C(D^!,\alpha^!)^{op}$
with  $C(D^!,\bar{\alpha}^!)$ in an obvious way,
we make of $\Cal{H}$ an $C(D,\alpha)$--$C(D^!,\bar{\alpha}^!)$
bimodule.
 
\smallskip

Assuming that we are in the situation of (3.3)
and taking the space $L_2(K_0)$ (rather than
$L_2(\widehat{K}_0)$) for the basic representation,
we will construct the Hermitean scalar products with the 
properties summarized in the Lemma 1.3.1. For further
details, see [Ri3]. 
We will assume that $K_0$ is a Lie group of the form
$\bold{R}^p\times \bold{Z}^q\times (\bold{R}/\bold{Z})^r\times F$
where $F$ is a finite group. Then one can define
the Schwartz space $S(K_0)$ consisting of
$C^{\infty}$--functions such that any polynomial
times any derivative of the function vanishes at infinity.
Rieffel's scalar products are first defined on
Schwartz's functions on $K_0$ and $\Cal{K}$ and then extended
to the appropriate completions. We will write
elements of $C(D,\alpha )$ as formal series
$F=\sum_{h\in D} a_he_{D,\alpha}(h)$ where $a_h$ are the
(non--commutative) Fourier coefficients defined
by $a_h=t(Fe(h)^*)$, $t$ is the normalized trace.
If Fourier coefficients form a Schwartz function
on $D$, $F$ will be called smooth. 

\smallskip

We start with the standard scalar
product on $L_2(K_0)$ (antilinear in the second argument)
which will be denoted $\langle *,*\rangle_{L_2}$
and put
for $\Phi,\Psi\in S(K_0)$:
$$
{}_D\langle \Phi ,\Psi\rangle :=\sum_{h\in D}
\langle \Phi ,e_{D,\alpha}(h)\Psi\rangle_{L_2}\,e_{D,\alpha}(h),
\eqno(3.12)
$$
$$
\langle \Phi ,\Psi\rangle_{D^!} :=\sum_{h\in D^!}
\langle e_{D^!,\alpha^!}(h)\Psi ,\Phi\rangle_{L_2}\,e_{D^!,\bar{\alpha}^!}(h).
\eqno(3.13)
$$
(Notice the appearance of both $\alpha^!$ and $\bar{\alpha}^!$
in the right hand side of (3.13)).

\smallskip

Before summarizing some results due to Rieffel, we
have to add a few words about the normalizations
of various Haar measures involved. Any Haar measure on $K_0$ will do;
on $\widehat{K}_0$ we take the respective Plancherel measure.
For the volumes of the respective fundamental
domains we will then have $|\Cal{K}/D|\cdot |\Cal{K}/D^!|=1.$

\medskip

\proclaim{\quad 3.4.1. Theorem} Denote by $M$
the completion of $S(K_0)$ with respect to
the operator norm
$\|{}_D\langle \Phi ,\Phi\rangle \|^{1/2}$.
Put 
$$
A=C(D,\alpha), \ B=C(D^!,\bar{\alpha}^!),\
{}_A\langle *,*\rangle =|\Cal{K}/D|\, {}_D\langle *,*\rangle,\
\langle *,*\rangle_B = \langle *,*\rangle_{D^!}.
$$
Then we have:

\smallskip

(a) $M$ is a finitely generated projective $A$--$B$ module
isomorphic to the range of a projection (both right and left).

\smallskip

(b) $A$ is the complete endomorphism ring of $M_B$.

\smallskip

(c) The scalar products defined above satisfy all
the identities (1.8)--(1.11).

\smallskip

(d) Let $t_B$ be the normalized trace on $B$
(zeroth Fourier coefficient). Then
$$
t_B([M_B])=|\Cal{K}/D^!|.
\eqno(3.14)
$$  
\endproclaim   

\smallskip

Notice that, contrary to the purely algebraic
context of Lemma 1.3.1 where (1.8)--(1.11) followed directly
from the definitions (1.6), (1.7), the deduction of (1.11)
from (3.12), (3.13) requires application of the Poisson summation formula.

\smallskip

For further details, see [Ri3], sections 2 and 3.

\medskip

{\bf 3.5. Vector Heisenberg group and classical theta functions.}
We return now temporarily to the setup of 3.1--3.2, involving
no additional lattice $D$ and explain the appearance
of the classical theta functions as matrix
coefficients of the basic Heisenberg representation.
We closely follow [Mu3], \S 3. 

\smallskip

We choose as $\Cal{K}$ the real vector space $V=\bold{R}^{2N}$.
Any element, say, $x\in V$ will be considered as a pair of
columns of height $N$: $x_1$ consisting of the first
$N$ coordinates of $x$ and $x_2$ consisting of the last
$N$ coordinates. Define the standard symplectic form on
$V$  and the cocycle $\psi$ (cf. (3.1)) by
$$
A(x,y)=x_1^ty_2-x_2^ty_1,\quad \psi (x,y)=e^{\pi iA(x,y)},
\eqno(3.15)
$$
so that
$$
\epsilon (x,y)=e^{2\pi iA(x,y)}.
\eqno(3.16)
$$
Having chosen a model $\Cal{H}$ of the basic representation
of the resulting Heisenberg group, Mumford defines
in $\Cal{H}$ a finite--dimensional family
of vectors $f_T\in \Cal{H}$ parametrized by the points
$T$ in the Siegel upper half space $\frak{H}_N$
consisting of complex symmetric $N\times N$
matrices with positive definite imaginary part.

\smallskip

In abstract terms, this is the space of all flat K\"ahler structures
on $\bold{R}^{2N}$ compatible with $A$. Such a structure
can be thought of, for example, as a pair
consisting of a complex structure $J$ and a positive
definite Hermitean form $H$ with the imaginary part $A$.

\smallskip

Any given $T$ determines directly the complex structure $J_T$:
it is given by thecomplex coordinates 
$\underline{x}_1,\dots ,\underline{x}_N$ on $V$:
$$
\underline{x}_i=\sum_j T_{ij}x_j^{(1)}+x_i^{(2)}
\eqno(3.17)
$$
where now  $x_j^{(1)}$ (resp. $x_i^{(2)}$) are the coordinates
of $x_1$ (resp. $x_2$). 
The values of the Hermitean form $H_T$ on the basic
vectors $e_j^{(2)}$ of the second half of $V$ are
$$
H_T(e_i^{(2)},e_j^{(2)})=(\roman{Im}\,T)_{ij}^{-1}.
\eqno(3.18)
$$

\smallskip

We can now define $f_T$ in the Mumford's first realization of the
fundamental representation:
$$
\Cal{H}:= L_2(\bold{R}^N),\ \quad (U_{(\lambda ,y_1,y_2)}f)(x)=
\lambda e^{2\pi ix_1^ty_2+\pi iy_1^tx_2} f(x_1+y_1)
\eqno(3.19)
$$
which is a specialization of (3.6) restricted to the subspace (3.8).
Namely, we have
$$
f_T(x)=e^{\pi i x^tTx}.
\eqno(3.20)
$$
The classical theta function is defined by
$$
\theta (\underline{x},T):= \sum_{n\in\bold{Z}^N}
e^{\pi in^tTn+2\pi in^t\underline{x}}.
\eqno(3.21)
$$
To express it as a matrix coefficient, Mumford 
introduces the distribution 
$$
e_{\bold{Z}}:=\sum_{n\in\bold{Z}^N} \delta_n
\eqno(3.22)
$$
and then checks that
$$
\langle U_{(1,x)}f_T,e_{\bold{Z}}\rangle =
c\,e^{\pi ix^t_1\underline{x}} \theta (\underline{x},T).
\eqno(3.23)
$$

(See [Ma3], Corollary 3.4).
\medskip

{\bf 3.6. Quantum theta functions.}
In this subsection I give a brief review of the formalism
of quantum theta functions introduced in [Ma1]
and further studied in [Ma2], [Ma3]. For details
and motivation, see [Ma3].
\smallskip

Consider the character group of a quantum torus $(D,\alpha )$. In this
subsection we will be interested in the space of
formal infinite linear combinations
of $e(h)=e_{D,\alpha}(h)$ which we will call formal functions.
Theta functions are defined as solutions
of functional
equations which can be invariantly described in terms of
another version
of Heisenberg group $\Cal{G}(D,\alpha )$ acting on this space: it consists
of all linear operators on formal functions
of the form
$$
\Phi\mapsto
c\,e(g)\,x^*(\Phi )\,e(h)^{-1} 
$$
where $c\in \bold{C}^*$, $g,h\in D$, 
$x\in T(D,1)(\bold{C}):=\roman{Hom}\,(D,\bold{C}^*)$ an arbitrary
point of the algebraic torus with the character group $D$,
$x^*$ is the shift automorphism multiplying
$e_{D,\alpha}(h)$ by $x(h)$.

\smallskip

Notice that such a shift $x^*$ generally does not respect the unitarity of
$e_{D,\alpha}(h)$ and cannot be extended to
the automorphisms of $C(D,\alpha )$ unless
the values of $x$ belong to $\bold{C}_1^*$.

\smallskip

We now define a (formal) theta multiplier 
for $(D,\alpha )$ as an injective homomorphism $\Cal{L}:\,
B\to \Cal{G} (D,\alpha )$ where $B$ is a free abelian group
of the same rank as $D$.

\smallskip

A quantum theta function with multiplier  $\Cal{L}$ is a formal
function on $T(D,\alpha )$ invariant with respect to the
action of (the image of) $B$.

\smallskip

$\Gamma (\Cal{L})$ is the linear space of theta functions with multiplier  $\Cal{L}$.

\smallskip

The theta functions constructed below will have coefficients
from the Schwartz space of $D$ and therefore will represent
smooth elements of $C(D,\alpha ).$ Their multipliers
will have the property $\roman{dim}\,\Gamma (\Cal{L})=1$.
In other words, as the classical $\theta (\underline{x},T)$,
our thetas will correspond only to the principal
polarizations. In order to get more general thetas one should consider 
more general Heisenberg groups in the sense of 3.1
into which $(D,\alpha )$ could be embedded.
 
\medskip

\proclaim{\quad 3.7. Theorem} We have
$$
{}_D\langle f_T,f_T\rangle = \frac{1}{\sqrt{2^N\roman{det}\, \roman{Im}\,T}}
\sum_{h\in D}e^{-\frac{\pi}{2}\underline{h}^t (\roman{Im}\,T)^{-1}
\underline{h}^*}e_{D,\alpha}(h),
\eqno(3.24)
$$
$$
\langle f_T,f_T\rangle_{D^!}= \frac{1}{\sqrt{2^N\roman{det}\, \roman{Im}\,T}}
\sum_{h\in D^!}e^{-\frac{\pi}{2}\underline{h}^t (\roman{Im}\,T)^{-1}
\underline{h}^*}e_{D^!,\bar{\alpha}^!}(h).
\eqno(3.25)
$$
Here $\underline{h} := Th_1+h_2$ (cf. (3.17)) and
$\underline{h}^* := \overline{T}h_1+h_2$.
These scalar products are quantum theta functions $\Theta_D,\Theta_{D^!}$
satisfying the following functional equations:
$$
\forall g\in D,\quad c_ge_{D,\alpha}(g)x^*_g(\Theta_D)=\Theta_D,
\eqno(3.26)
$$
$$
\forall g\in D^!,\quad c_g^!e_{D^!,\bar{\alpha}^!}(g)x^{!*}_g(\Theta_{D^!})=
\Theta_{D^!}
\eqno(3.27)
$$
where
$$
c_g=e^{\frac{3\pi}{2}\underline{g}^t (\roman{Im}\,T)^{-1}
\underline{g}^*},\quad x^*_g(e_{D,\alpha}(h))=e^{X_g(h)}e_{D,\alpha}(h),
$$
$$
X_g(h)=-\pi\,\roman{Re}\,(\underline{g}^t (\roman{Im}\,T)^{-1}
\underline{h}^*) -\pi i\,A(g,h),
\eqno(3.28)
$$
$$
c_g^!=e^{\frac{3\pi}{2}\underline{g}^t (\roman{Im}\,T)^{-1}
\underline{g}^*},\quad
x^{!*}_g(e_{D^!,\bar{\alpha}}(h))=e^{X_g^!(h)}e_{D^!,\bar{\alpha}}(h), 
$$
$$
X_g^!(h)=-\pi\,\roman{Re}\,(\underline{g}^t (\roman{Im}\,T)^{-1}
\underline{h}^*) +\pi i\,A(g,h).
\eqno(3.29)
$$

\endproclaim

\smallskip

{\bf Proof.} We will check (3.24) and (3.26); the other two
formulas can be treated similarly.

\smallskip

The general formula (3.12) must be specialized to our case
$K_0=\bold{R}^N,$ the first half of $\Cal{K}.$ For the
$L_2$--scalar product we take $\int \Phi\overline{\Psi}dx_1$
where $dx_1$ is the standard Haar measure. From (3.10) and
(3.19) it follows that
$$
(e_{D,\alpha}(h)\Psi )(x_1)= 
e^{2\pi ix_1^th_2+\pi ih_1^th_2} \Psi (x_1+h_1).
$$
Hence
$$
{}_D\langle \Phi ,\Psi\rangle =\sum_{h\in D}
e^{-\pi ih_1^th_2}\int \Phi (x_1)\overline{\Psi (x_1+h_1)}
e^{-2\pi ix_1^th_2}dx_1\cdot e_{D,\alpha}(h).
$$
Putting here $\Phi =\Psi =f_T$ (see (3.20)), we get:
$$
{}_D\langle f_T ,f_T\rangle =\sum_{h\in D}
e^{-\pi ih_1^th_2}\int 
e^{\pi i[x_1^tTx_1-(x_1^t+h_1^t)\overline{T}(x_1+h_1)-2x_1^th_2]}
dx_1\cdot e_{D,\alpha}(h).
\eqno(3.30)
$$
The exponential expression under the integral sign in (3.30)
can be represented as $e^{-(q(x_1)+l_h(x_1)+c_h)}$ where
$$
q(x_1)=2\pi\,x_1^t\,(\roman{Im}\,T)^{-1}\,x_1,\ 
l_h(x_1)=2\pi i\,(h_1^t\overline{T}+h_2^t)\,x_1,\
c_h=\pi ih_1^t\overline{T} h_1 .
\eqno(3.31)
$$
Denote
$$
\lambda_h=\frac{i}{2\,\roman{Im}\,T}[\overline{T}h_1+h_2].
\eqno(3.32)
$$
Then we have
$$
q(x_1+\lambda_h)-q(\lambda_h) =q(x_1)+l_h(x_1)
$$
and therefore
$$
\int e^{-(q(x_1)+l_h(x_1)+c_h)}dx_1=
e^{-c_h+q(\lambda_h)}\int e^{-q(x_1+\lambda_h)} dx_1=
e^{-c_h+q(\lambda_h)}\frac{\pi^{N/2}}{\sqrt{\roman{det}\,q }}.
\eqno(3.33)
$$
Putting (3.30)--(3.33) together, we get (3.24).

\smallskip

The equation (3.26) is checked by a straightforward
computation: putting $Q(h)=\dfrac{\pi}{2}\,h_1^t\,(\roman{Im}\,T)^{-1}\,h_1$
we have
$$
c_ge_{D,\alpha}(g)x_g^*(\Theta_D)=c_g\sum_{h\in D}
e^{-Q(h)+X_g(h)+\pi i\,A(g,h)}e_{D,\alpha}(g+h)=
$$
$$
c_g\sum_{h\in D}
e^{-Q(h-g)+X_g(h-g)+\pi i\,A(g,h-g)}e_{D,\alpha}(h)=
$$
$$
c_ge^{-Q(g)+X_g(-g)}\sum_{h\in D}
e^{-Q(h)+Q(g)+X_g(h)+\pi i\,A(g,h)}e_{D,\alpha}(h)= \Theta_D.
$$

\newpage

\centerline{\bf Bibliography}

\medskip

[ArLR] J.~Arledge, M.~Laca, I.~Raeburn. {\it Crossed products by semigroups
of endomorphisms and the Toeplitz algebras of ordered groups.}
Documenta Math., 2 (1997), 115--138.

\smallskip

[AsSch] A.~Astashkevich, A.~S.~Schwarz. {\it Projective modules over
non--commutative tori: classification of modules with
constant curvature connection.} e--Print QA/9904139.

\smallskip

[BG] V.~Baranovsky, V.~Ginzburg. {\it Conjugacy classes
in loop groups and $G$--bundles on elliptic curves.}
Int. Math. Res. Notices, 15 (1996), 733--751.

\smallskip

[BEG] V.~Baranovsky, S.~Evens, V.~Ginzburg. {\it Representations
of quantum tori and double--affine Hecke algebras.} e--Print
math.RT/0005024

\smallskip

[Bl] B.~Blackadar. {\it $K$--theory for operator algebras.}
Springer Verlag, 1986.

\smallskip

[Bo1] F.~Boca. {\it The structure of higher--dimensional
noncommutative tori and metric diophantine approximation.}
J. Reine Angew. Math. 492 (1997), 179--219.

\smallskip

[Bo2] F.~Boca. {\it Projections in rotation algebras and theta functions.}
Comm. Math. Phys., 202 (1999), 325--357.

\smallskip 

[BoCo] J.~Bost, A.~Connes. {\it Hecke algebras, type III factors and phase
transitions with spontaneous symmetry breaking in number theory.}
Selecta Math., 3 (1995), 411--457.

\smallskip

[BrGrRi] L.~Brown, P.~Green, M.~A.~Rieffel.
{\it Stable isomorphism and strong Morita equivalence
of $C^*$--algebras.} Pacific J.~Math., 71 (1977),
349--363.

\smallskip

[Coh1] P.~Cohen. {\it A $C^*$--dynamical system with Dedekind zeta partition function
and spontaneous symmetry breaking.}
J. de Th. de Nombres de Bordeaux, 11 (1999), 15--30.

\smallskip

[Coh2] P.~Cohen. {\it Quantum statistical mechanics and number theory.}
In: Algebraic Geometry: Hirzebruch 70. Contemp. Math.,
vol 241,  AMS, Pridence RA, 1999, 121--128.

\smallskip

[Co1] A.~Connes. {\it Noncommutative geometry.} Academic Press,
1994.

\smallskip

[Co2] A.~Connes. {\it $C^*$--algebras et g\'eom\'etrie
diff\'erentielle.} C.~R.~Ac.~Sci.~Paris, t.~290 (1980), 599--604.

\smallskip

[Co3] A.~Connes. {\it Trace formula in noncommutative geometry
and the zeros of the Riemann zeta function.} Selecta Math., New Ser., 5 (1999),
29--106.

\smallskip

[Co4] A.~Connes. {\it Noncommutative geometry and the Riemann
zeta function.} In: Mathematics: Frontiers and Perspectives,
ed. by V.~Arnold et al., AMS, 2000, 35--54.

\smallskip

[Co5] A.~Connes. {\it Noncommutative Geometry Year 2000.}
e--Print   math.QA/0011193

\smallskip

[CoDSch] A.~Connes, M.~Douglas, A.~Schwarz. {\it Noncommutative
geometry and Matrix theory: compactification on tori.}
Journ. of High Energy Physics, 2 (1998).

\smallskip
[Dav] K.~R.~Davidson. {\it $C^*$--algebras by example.} AMS, 1996.

\smallskip

[Dar] H.~Darmon. {\it Stark--Heegner points over real quadratic fields.}
Contemp. Math., 210 (1998), 41--69.

\smallskip

[DE] {\it The Grothendieck Theory of Dessins d'Enfants
(Ed. by L.~Schneps).} London MS Lect. Note series, 200,
Cambridge University Press, 1994.

\smallskip

{DiSch] M.~Dieng, A.~Schwarz. {\it Differential
and
complex geometry of two--dimensional noncommutative tori.}
Preprint, 2002.

\smallskip
[Dr] V.~G.~Drinfeld. {\it On quasi--triangular quasi--Hopf
algebras and some groups closely associated
with $Gal (\overline{\bold{Q}}/\bold{Q})$.}
Algebra and Analysis 2:4 (1990); Leningrad Math. J.
2:4 (1991), 829--860.

\smallskip

[El1] G.~A.~Elliott. {\it On the $K$--theory of the $C^*$--algebra
generated by a projective representation of a torsion--free discrete
abelian group.} In: Operator Algebras
and Group Representations I, Pitman, London, 1984, 157--184.

\smallskip

[El2] G.~A.~Elliott. {\it The diffeomorphism groups of the irrational
rotation $C^*$--algebra.} C.~R.~Math.~Rep.~Ac.~Sci.~Canada,
8 (1986), 329--334.

\smallskip

[ElEv] G.~A.~Elliott, D.~E.~Evans. {\it The structure of the
irrational rotation $C^*$--algebra.} Ann. Math., 138 (1993), 477--501.

\smallskip

[FK] L.~D.~Faddeev, R.~M.~Kashaev. {\it Quantum dilogarithm}.
e--Print hep-th/9310070.

\smallskip

[Ga] C.~F.~Gauss.  {\it Zur Kreistheilung.} In: Werke, Band 10, p.4.
Georg Olms Verlag, Hildsheim--New York, 1981.

\smallskip

[GoHaJo] F.~Goodman, P.~de la Harpe, V.~Jones. {\it Coxeter graphs and
towers of algebras.} Springer, 1989.

\smallskip

[HaL] D.~Harari, E.~Leichtnam. {\it Extension du
ph\'enom\`ene de brisure spontan\'ee de sym\'etrie
de Bost--Connes au cas des corps globaux
quelconques.} Selecta Math., 3 (1997), 205--243.

\smallskip

[He1] E.~Hecke. {\it \"Uber die Kroneckersche Grenzformel
f\"ur reelle quadratische K\"orper und die Klassenzahl
relativ--abelscher K\"orper.} Verhandl. d. 
Naturforschender Gesell. i. Basel, 28 (1917),
363--373. (Math. Werke, pp. 198--207, Vandenberg \& Ruprecht,
G\"ottingen, 1970).

\smallskip

[He2] E.~Hecke. {\it Zur Theorie der elliptischen  
Modulfunktionen.} Math.~Annalen, 97 (1926), 210--243.
(Math. Werke, pp. 428--460, Vandenberg \& Ruprecht,
G\"ottingen, 1970).

\smallskip

[Her] G.~Herglotz. {\it \"Uber die Kroneckersche Grenzformel
f\"ur reelle quadratische K\"orper I, II.} Berichte
\"uber die Verhandl. S\"achsischen Akad. der Wiss.
zu Leipzig, 75 (1923), 3--14, 31--37.

\smallskip

[Jo1] V.~Jones. {\it Index for subfactors.} Inv. Math., 72:1 (1983), 1--25.

\smallskip

[Jo2] V.~Jones. {\it Index for subrings of rings.}
Contemp. Math. 43, AMS(1985), 181--190.

\smallskip

[Jo3] V.~Jones. {\it Fusion en alg\`ebres de von Neumann et groupes de
lacets (d'apr\`es A.~Wassermann)}. S\'eminaire Bourbaki, no. 800
(Juin 1995), 20 pp.

\smallskip

[JoSu] V.~Jones, V.~Sunder. {\it Introduction to subfactors.}
Lond. Math. Soc. Lecture Note Series, 234 (1997).

\smallskip

[KaNi1] L.~Kadison, D.~Nikshych. {\it Frobenius extensions
and weak Hopf algebras.} e--Print math/0102010 .

\smallskip

[KaNi2] L.~Kadison, D.~Nikshych. {\it Hopf algebra actions
on strongly separable extensions of depth two.}
e--Print math/0107064 .

\smallskip

[KoSo] M.~Kontsevich, Y.~Soibelman. {\it Homological
mirror symmetry and torus fibrations.} e--Print  
math.SG/0011041

\smallskip

[La] S.~Lang. {\it Elliptic Functions.} Addison--Wesley, 1973.

\smallskip

[Lin] Qing Lin. {\it Cut--down method in the inductive limit
decomposition of non--commutative tori, III: A complete
answer in 3--dimension.} Comm.~Math.~Phys., 179 (1996), 555--575.

\smallskip

[LoS] P.~Lochak, L.~Schneps. {\it A cohomological interpretation
of the Grothendieck--Teichm\"uller group.} Inv. Math.,
127 (1997), 571--600.

\smallskip

[Ma1] Yu.~Manin. {\it  Quantized theta--functions.} In: Common
Trends in Mathematics and Quantum Field Theories (Kyoto, 1990), 
Progress of Theor. Phys. Supplement, 102 (1990), 219--228.

\smallskip

[Ma2] Yu.~Manin. {\it Mirror symmetry and quantization of abelian varieties.}
In: Moduli of Abelian Varieties, ed. by C.~Faber et al.,
Progress in Math., vol. 195, Birkh\"auser, 2001, 231--254.
e--Print  math.AG/0005143

\smallskip

[Ma3] Yu.~Manin. {\it Theta functions, quantum tori and Heisenberg groups}.
Lett. in Math. Physics, 56 (2001), 295--320.
e--Print math.AG/001119

\smallskip

[Ma4] Yu.~Manin. {\it Von Zahlen und Figuren.} 27 pp. e--Print math.AG/0201005

\smallskip

[MaMar] Yu.~Manin, M.~Marcolli. {\it Continued fractions, modular symbols, and non-commutative geometry.} e--Print math.NT/0102006

\smallskip

[Mu1] D.~Mumford. {\it On the equations defining abelian
varieties I.} Inv. Math. 1 (1966), 287 --354.

\smallskip

[Mu2] D.~Mumford. {\it An algebro--geometric construction
of commuting operators and of solutions to the Toda
lattice equation, Korteweg deVries equation and related
non--linear equations.} In: Proc. of Int. Symp. on
Algebraic Geometry, Kyoto, 1977, 115--153.

\smallskip

[Mu3] D.~Mumford (with M.~Nori and P.~Norman). 
{\it Tata Lectures on Theta III.} Progress in Math., vol.~97,
Birkh\"auser, 1991.

\smallskip

[NiVa] D.~Nikshych, L.~Vainerman. {\it Finite quantum
groupoids and their applications.} e--Print math.QA/0006057.

\smallskip

[PimV] M.~Pimsner, D.~Voiculescu. {\it Exact sequences for $K$--groups
and $EXT$--groups of certain cross--product $C^*$--algebras.}
J.~Operator Theory, 4 (1980), 93--118.

\smallskip

[Po1] A.~Polishchuk. {\it Indefinite theta series
of signature (1,1) from the point of view of
homological mirror symmetry.} e--Print 
math.AG/0003076 .

\smallskip

[Po2] A.~Polishchuk. {\it A new look at Hecke's indefinite theta series.}
e--Print  math.AG/0012005 .

\smallskip

[RaW] I.~Raeburn, D.~Williams. {\it Morita equivalence and
continuous--trace $C^*$--algebras.} Math. Surveys and Monographs,
vol. 60, AMS, 1998.

\smallskip

[Ri1] M.~A.~Rieffel. {\it Strong Morita equivalence of
certain transformation group $C^*$--algebras.}
Math. Annalen, 222 (1976), 7--23.

\smallskip

[Ri2] M.~A.~Rieffel. {\it Von Neumann algebras associated with pairs of
lattices in Lie groups.} Math.~Ann., 257 (1981), 403--418.

\smallskip

[Ri3] M.~A.~Rieffel. {\it $C^*$--algebras associated with irrational rotations.}
Pacific J.~Math., 93 (1981), 415--429.

\smallskip

[Ri4] M.~A.~Rieffel. {\it The cancellation theorem for projective
modules over irrational rotation $C^*$--algebras.}
Proc.~Lond.~Math.~Soc. (3), 47 (1983), 285--303.

\smallskip

[Ri5] M.~A.~Rieffel. {\it Projective modules over higher--dimensional
non--commutative tori.} Can.~J.~Math., vol.~XL, No.~2 (1988), 257--338.

\smallskip

[Ri6] M.~A.~Rieffel. {\it Non--commutative tori --- a case
study of non--commutative differential manifolds.}
In: Cont.~Math., 105 (1990), 191--211.

\smallskip

[RiSch] M.~A.~Rieffel, A.~Schwarz. {\it Morita equivalence
of multidimensional non--commutative tori.} 
Int. J. Math., 10 (1999), 289--299. e--Print   math.QA/9803057

\smallskip

[Ro1] A.~Rosenberg. {\it Non--commutative algebraic geometry
and representations of quantized algebras.} Kluwer
Academic Publishers, 1995.

\smallskip

[Ro2] A.~Rosenberg. {\it Non--commutative schemes.}
Composition Math., 112 (1998), 93--125.

\smallskip

[Sch1] A.~Schwarz. {\it Morita equivalence and duality.}
Nucl.~Phys., B 534 (1998), 720--738.

\smallskip

[Sch2] A.~Schwarz. {\it Theta--functions on non--commutative tori.}
e--Print math/0107186

\smallskip

[Se] J.~P.~Serre. {\it Complex Multiplication.} In: 
Algebraic Number Fields, ed. by J.~Cassels, A.~Fr\"olich.
Academic Press, NY 1977, 293--296.

\smallskip

[So] Y.~Soibelman. {\it Quantum tori, mirror symmetry
and deformation theory.} Lett. in Math. Physics,
56 (2001), 99--125. e--Print math.QA/0011162.

\smallskip

[St1] H.~M.~Stark. {\it $L$--functions at $s=1$. III. Totally real
fields and Hilbert's Twelfth Problem.} Adv. Math., 22 (1976), 64--84.

\smallskip
[St2] H.~M.~Stark. {\it $L$--functions at $s=1$. IV. First
derivatives at $s=0$.} Adv. Math., 35 (1980), 197--235.

\smallskip

[Ste] P.~Stevenhagen. {\it Hilbert's 12th problem, Complex
Multiplication and Shimura reciprocity.} In: Class Field Theory --
Its Centenary and Prospect. Adv. Studies in Pure Math.,
30 (2001), 161--176.

\smallskip

[Ta] J.~Tate. {\it Les conjectures de Stark sur les fonctions
$L$ d'Artin en $s=0$.} Progress in Math., vol. 47, Birkh\"auser, 1984.

\smallskip

[U] A.~Unterberger. {\it Quantization and non--holomorphic
modular forms.} Springer LNM, 1742 (2001).

\smallskip

[Wa] Y.~Watatani. {\it Index for $C^*$--subalgebras.}
Mem. AMS, vol. 83, Nr. 424, Providence, RA, 1990.

\smallskip

[Wi] E.~Witten. {\it Overview of $K$--theory applied to strings.}
e--Print hep-th/0007175

\smallskip

[Za1] D.~Zagier. {\it A Kronecker limit formula for real quadratic field.}
Math. Ann., 213 (1975), 153--184.

\smallskip

[Za2] D.~Zagier. {\it Valeurs des fonctions zeta des corps
quadratiques r\'eels aux entiers n\'egatifs.} Ast\'erisque
41--42 (1977), 135--151.

\smallskip

[Zap] L.~Zapponi. {\it Dessins d'enfants en genre 1.}
In: Geometric Galois Actions, ed. L.~Schneps, P.~Lochak.

\bigskip

e-mail: manin\@mpim-bonn.mpg.de

\enddocument